\documentclass[10pt,a4paper,oneside]{amsart}
\usepackage{amsfonts, amscd, amsmath, amssymb}
\usepackage[dvips]{graphicx}
\newcommand\globcnt{subsection}
\newcommand\globsubcnt{subsubsection}

\theoremstyle{plain}
\newtheorem*{ThStabA}{Theorem~\ref{th:Stab}}

\newtheorem{thm}[\globcnt]{Theorem}
\newtheorem{lem}[\globsubcnt]{Lemma}
\newtheorem{prop}[\globsubcnt]{Proposition}
\newtheorem{cor}[\globsubcnt]{Corollary}

\newtheorem{claim}[\globsubcnt]{Claim}
\theoremstyle{definition}
\newtheorem{defn}[\globsubcnt]{Definition}
\newtheorem{exmp}[\globsubcnt]{Example}
\newtheorem{rem}[\globsubcnt]{Remark}

\makeatletter
\@addtoreset{subsubsection}{section}
\@addtoreset{subsubsection}{subsection}
\@addtoreset{equation}{section}
\makeatother

\providecommand{\eqref}[1]{(\ref{#1})}

\newcommand\CCC{{\mathbb C}}

\newcommand\RRR{{\mathbb R}}

\newcommand\ZZZ{{\mathbb Z}}

\newcommand\FFF{{\mathcal F}}

\newcommand\VV{{\mathcal V}}

\newcommand{\id}{\mathrm{id}}

\newcommand{\supp}{\mathrm{supp\,}}

\newcommand\Orbit{\mathcal{O}}
\newcommand\Stab{\mathcal{S}}
\newcommand{\Diff}{\mathcal{D}}

\newcommand\pnt{z}
\newcommand\xpnt{x}

\newcommand\eps{\varepsilon}

\newcommand\manif{M}
\newcommand\nmanif{N}
\newcommand\bfunc{\beta}
\newcommand\cfunc{\gamma}
\newcommand\dfunc{\delta}
\newcommand\ofunc{\omega}
\newcommand\nfunc{\nu}

\newcommand\homf{f}

\newcommand\difR{\phi}
\newcommand\difP{\difR}

\newcommand\nbh{U}

\newcommand\mrsfunc{f}

\newcommand\dimM{m}

\newcommand\difM{h}

\newcommand\singf{\Sigma_\mrsfunc}

\newcommand\Psp{P} 
\newcommand\Circle{S^1}

\newcommand\afunc{\alpha}
\newcommand\bafunc{\bar\alpha}

\newcommand\bdifM{\bar{\difM}}

\newcommand\sct{s}

\newcommand\mnfunc{u}
\newcommand\dmnfunc{v}
\newcommand\tmnfunc{\omega}

\newcommand\qfunc{q}
\newcommand\gfunc{g}
\newcommand\hfunc{h} 
\newcommand\crvmap{k}

\newcommand\anbh{V}
\newcommand\bdifR{\bar\difR}

\newcommand\barx{\overline{x}}

\newcommand\cnv{i}

\newcommand\germm{C^{\infty}_{\pnt}(\manif)}
\newcommand\JFP[2]{\Delta(#1,#2)}
\newcommand\JFzero[1]{\JFP{#1}{0}}

\newcommand\Jacidfzero{\JFzero{\mrsfunc}}
\newcommand\Jacidgzero{\JFzero{\gfunc}}

\newcommand\Jacidfz{\JFP{\mrsfunc}{\pnt}}

\newcommand\prMPP{p}

\newcommand\sectMPP{\Theta}
\newcommand\difShift{\theta}
\newcommand\homood{\alpha}
\newcommand\homoodinv{\beta}

\newcommand\dntwo{\Delta^{\nv-2}}
\newcommand\dnone{\Delta^{\nv-1}}
\newcommand\dn{\Delta^{\nv}}
\newcommand\Pcr{E}
\newcommand\pcr{e}
\newcommand\tpcr{\widetilde{\pcr}}
\newcommand\padj{q}
\newcommand\adjd{c}
\newcommand\cfuncinv{\cfunc^{-1}}
\newcommand\hconf{\xi}

\newcommand\Nbh{\mathcal{\nbh}}

\newcommand\dfval[2]{#2.#1}
\newcommand\dfvalx[3]{\dfval{#1}{#2}(#3)}

\newcommand\sfunc{\sigma}

\newcommand\rnbh{V}
\newcommand\mnbh{U}

\newcommand\vectformat[1]{#1}
\newcommand\flowformat[1]{\mathcal{#1}}

\newcommand\Fld{\vectformat{H}}

\newcommand\fFld{\vectformat{F}}
\newcommand\gFld{\vectformat{G}}
\newcommand\hFld{\vectformat{H}}

\newcommand\rFld{\vectformat{F}}
\newcommand\mFld{\vectformat{H}}

\newcommand\rflow{\flowformat{F}}
\newcommand\mflow{\flowformat{H}}

\newcommand\caseA{{\rm a)}}
\newcommand\caseB{{\rm b)}}

\newcommand\indx{k}

\newcommand\DiffP{\Diff_{\Psp}}
\newcommand\DiffR{\Diff_{\RRR}}

\newcommand\DiffM{\Diff_{\manif}}

\newcommand\DiffSpl{\Diff_{\Circle}^{+}}

\newcommand\DiffSe{\Diff_{\Circle}^{e}}

\newcommand\DiffPi{\Diff_{\Psp}'}

\newcommand\DiffRf{\Diff_{\RRR}^{[1,n]}}
\newcommand\DiffSf{\Diff_{\Circle}^{+}}

\newcommand\DiffRfe{\Diff_{\RRR}^{e}}
\newcommand\DiffSfe{\Diff_{\Circle}^{e}}
\newcommand\DiffSfE{\Diff_{\Circle}^{E}}
\newcommand\DiffPfe{\Diff_{\Psp}^{e}}

\newcommand\DiffMRf{\Diff_{\manif\RRR}}
\newcommand\DiffMSf{\Diff_{\manif\Circle}}
\newcommand\DiffMPf{\Diff_{\manif\Psp}}

\newcommand\OrbfMRi{\Orbit_{\manif\RRR}'(\mrsfunc)}

\newcommand\OrbfMSi{\Orbit_{\manif\Circle}'(\mrsfunc)}
\newcommand\StabfMPi{\Stab_{\manif\Psp}'(\mrsfunc)}
\newcommand\OrbfMPi{\Orbit_{\manif\Psp}'(\mrsfunc)}

\newcommand\StabfMP{\Stab_{\manif\Psp}(\mrsfunc)}
\newcommand\StabfMR{\Stab_{\manif\RRR}(\mrsfunc)}
\newcommand\StabfMS{\Stab_{\manif\Circle}(\mrsfunc)}
\newcommand\StabfM{\Stab_{\manif}(\homf)}

\newcommand\OrbfMP{\Orbit_{\manif\Psp}(\mrsfunc)}
\newcommand\OrbfMR{\Orbit_{\manif\RRR}(\mrsfunc)}
\newcommand\OrbfMS{\Orbit_{\manif\Circle}(\mrsfunc)}
\newcommand\OrbfM{\Orbit_{\manif}(\homf)}

\newcommand\excf{E_{\mrsfunc}}

\newcommand\smr{C^{\infty}(\manif,\RRR)}
\newcommand\sms{C^{\infty}(\manif,\Circle)}
\newcommand\smone{C^{\infty}(\manif,\Psp)}

\newcommand\pStabfMS{\widetilde{\Stab}_{\manif\Circle}(\mrsfunc)}
\newcommand\kcnv{c}
\newcommand\dcnv{d}

\newcommand\Rm{\RRR^{\dimM}}
\newcommand\Rn{\RRR^{\nv}}

\newcommand\minid[1]{\mathfrak{m}(#1)}
\newcommand\mmm{\minid{\Rm}}
\newcommand\mmmRm{\minid{\Rm}}

\newcommand\mmmR{\minid{\RRR}}

\newcommand\smzr{C^{\infty}_{0}(\RRR)}
\newcommand\smzrm{C^{\infty}_{0}(\Rm)}
\newcommand\Vdf[1]{L({#1})}
\newcommand\hmV[1]{\tau_{#1}}

\newcommand\Interv{I}

\newcommand\mfunc{\mu}
\newcommand\vars{s}

\newcommand\nv{n}
\newcommand\kv{k}

\newcommand\lnbh[1]{V_{#1}}
\newcommand\lFld[1]{\vectformat{H}_{#1}}
\newcommand\lflow[1]{\flowformat{H}_{#1}}
\newcommand\lfunc[1]{\sigma_{#1}}
\newcommand\ldif[1]{h_{#1}}
  
\newcommand\elnbh[1]{U_{#1}} 
\newcommand\elnbhp[1]{U^{+}_{#1}} 
\newcommand\elnbhm[1]{U^{-}_{#1}} 
\newcommand\elFld[1]{\vectformat{G}_{#1}}
\newcommand\elflow[1]{\flowformat{G}_{#1}}
\newcommand\elfunc[1]{\rho_{#1}}
\newcommand\eldif[1]{g_{#1}}

\newcommand\exlev[1]{L_{#1}} 
\newcommand\exlevel{L}

\newcommand\jj{j}
\newcommand\elFldref[1]{\vectformat{G}'_{#1}}
\newcommand\elnbhref[1]{U'_{#1}}

\newcommand\condBndCrVal{{\rm(V)}}
\newcommand\condJacIda[1]{{\rm J(${#1}$)}} 
\newcommand\condJacId{{\rm (J)}} 
\newcommand\LT{{\rm L}}

\newcommand\ConfSp{\FFF_{\nv}(\Circle)}
\newcommand\Tn{T^{\nv}}
\newcommand\impcr{F_{\nv}(\Circle)}

\newcommand\acnv{a}

\newcommand\cnbh{W}
\newcommand\crset{C}
\newcommand\aNbh{\VV}

\newcommand\genzn{\nu}

\newcommand\Zk{\ZZZ_{\kcnv}}
\newcommand\traj[1]{\omega_{#1}}

\newcommand\pfact{\tau}
\newcommand\pOrbfMS{\widetilde{\Orbit}_{\manif\Circle}(\mrsfunc)}
\newcommand\prone{p_1}
\newcommand\crvmaptwo{\crvmap_2}
\newcommand\bij{\mu}

\newcommand\Gfunc{G}
\newcommand\factn{\nu}

\newcommand\DiffRmz{\Diff_{0}(\Rm)}
\newcommand\DiffRz{\Diff_{0}(\RRR)}
\newcommand\StabfRmR{\Stab_{\Rm,\RRR}(\mrsfunc)}

\newcommand\adif{\phi}
\newcommand\bdif{\psi}
\newcommand\cdif{\xi}

\newcommand\dgfunc[1]{\beta_{#1}}
\newcommand\shfunc[1]{\sigma_{\!#1}}

\newcommand\adgfn{\dgfunc{\adif}}

\newcommand\ashfn{\shfunc{\adif}}
\newcommand\bshfn{\shfunc{\bdif}}
\newcommand\cshfn{\shfunc{\cdif}}

\newcommand\mideal[1]{I_{#1}}
\begin{document}
\author{Sergey Maksymenko}
\title{Stabilizers and orbits of smooth functions}

\address{
 Topology Department,
 Institute of Mathematics, NAS of Ukraine,
 Tereshchenkivska str. 3, 01601 Kyiv, Ukraine,
 e-mail:\texttt{maks@imath.kiev.ua},
 phone: (+380 44) 2345150,
 fax: (+380 44) 2352010 
}

\maketitle

{
\footnotesize 
{\bf Abstract.}
Let $\mrsfunc:\Rm\to\RRR$ be a smooth function such that $\mrsfunc(0)=0$.
We give a condition \condJacIda{\id}\ on $\mrsfunc$ when for arbitrary preserving orientation diffeomorphism $\difR:\RRR\to\RRR$ such that $\difR(0)=0$ the function $\difR\circ\mrsfunc$ is right equivalent to $\mrsfunc$, i.e.\!\! there exists a diffeomorphism $\difM:\Rm\to\Rm$ such that $\difR\circ\mrsfunc=\mrsfunc\circ\difM$ at $0\in\Rm$.
The requirement is that $\mrsfunc$ belongs to its Jacobi ideal.
This property is rather general: it is invariant with respect to the stable equivalence of singularities, and holds for non-degenerated, simple, and many other singularities.

We also globalize this result as follows.
Let $\manif$ be a smooth compact manifold, $\mrsfunc:\manif\to[0,1]$ a surjective smooth function, $\DiffM$ the group of diffeomorphisms of $\manif$, and $\DiffR^{[0,1]}$ the group of diffeomorphisms of $\RRR$ that have compact support and leave $[0,1]$ invariant.
There are two natural right and left-right actions of $\DiffM$ and $\DiffM\times\DiffR^{[0,1]}$ on $\smr$.
Let $\StabfM$, $\StabfMR$, $\OrbfM$, and $\OrbfMR$ be the corresponding stabilizers and orbits of $\mrsfunc$ with respect to these actions.
We prove that if $\mrsfunc$ satisfies \condJacIda{\id}\ at each critical point and has additional mild properties, then the following homotopy equivalences hold: $\StabfM\approx\StabfMR$ and $\OrbfM\approx\OrbfMR$.
Similar results are obtained for smooth mappings $\manif\to\Circle$.

\medskip

{\bf R\'esum\'e.}
Soit $\mrsfunc:\Rm\to\RRR$ une application diff\'erentiable telle que $\mrsfunc(0)=0$.
On introduit une condition \condJacIda{\id}\ et preuve que si $\mrsfunc$ lui satisfait, alors pour chaque diff\'eomorphisme $\difR:\RRR\to\RRR$ tel que $\difR(0)=0$ et $\difR$ pr\'eserve l'orientation de $\RRR$ il existe un diff\'eomorphisme $\difM:\Rm\to\Rm$ tel que $\difR\circ\mrsfunc=\mrsfunc\circ\difM$ en $0\in\Rm$.
La condition \condJacIda{\id}\ requ\'erit que $\mrsfunc$ appartient \`a son id\'eal Jacobien en $0$.
Cette propri\'et\'e est invariante par rapport \`a l'equivalance stable des singularit\'es. 
Aussi les singularit\'es non-degenerates et simples satisfont \condJacIda{\id}. 

Ce r\'esultat local implique le th\'eor\`eme global suivant.
Soient $\manif$ une vari\'et\'e diff\'eren\-tiable, $\mrsfunc:\manif\to[0,1]$ une application diff\'erentiable surjective, 
$\DiffM$ le groupe des diff\'eomorphismes de $\manif$ et $\DiffR^{[0,1]}$ le groupe des diff\'eomorphismes de $\RRR$ de support compact et pr\'eservient $[0,1]$. 
Il y a deux actions droite et gauche-droite de $\DiffM$ et $\DiffM\times\DiffR^{[0,1]}$ respectivement sur $\smr$.
Soient $\StabfM$, $\StabfMR$, $\OrbfM$ et $\OrbfMR$ les stabilisateurs et orbites correspondantes de $\mrsfunc$ par rapport ces actions.
On preuve que si $\mrsfunc$ satisfait \condJacIda{\id}\ en chaque point critique et aussi une outre condition naturelle, alors on a les \'equivalences homotopies suivantes: $\StabfM\approx\StabfMR$ et  
$\OrbfM\approx\OrbfMR$.
Les r\'esultats ressemblantes sont obtenus pour les applications diff\'erentiables $\manif\to\Circle$.
}



\section{Introduction}
Let $\manif$ be a smooth ($C^{\infty}$) connected compact $\dimM$-dimensional manifold,
$\Psp$ either the real line $\RRR$ or the circle $\Circle$, and $\DiffM$ and $\DiffP$ the groups of diffeomorphisms of $\manif$ and $\Psp$ respectively.
There are natural {\em left\/} actions of the groups $\DiffM$ and
$\DiffM \times \DiffP$ on $\smone$ defined by the following formulas:
if $\mrsfunc\in\smone$,
$\difM\in\DiffM$, and $\difP\in\DiffP$, then
\begin{gather}
 \difM \cdot \mrsfunc  = \mrsfunc \circ \difM^{-1}
    \label{equ:act_M}\\
 (\difM,\difP) \cdot \mrsfunc  = \difP \circ \mrsfunc \circ \difM^{-1}.
    \label{equ:act_MP}
\end{gather}
These actions are often called \emph{right\/} and \emph{left-right} respectively. 
They were studied by many authors, see e.g.~\cite{AVG,LevineThom} for references.

For $\mrsfunc\in\smone$ let 
$\StabfM=\{\difM \in \DiffM \ | \ \mrsfunc=\mrsfunc\circ\difM \}$ be the stabilizer and 
$\OrbfM=\{\mrsfunc\circ\difM^{-1} \ | \ \difM\in\DiffM\}$ the orbit of $\mrsfunc$ under the right action~\eqref{equ:act_M}.

If $\DiffPi$ is a subgroup of $\DiffP$, then the left-right action of the group $\DiffM\times\DiffPi$ on $\smone$ is well-defined by~\eqref{equ:act_MP}.

Let also 
$\StabfMPi=\{(\difM,\difP) \in \DiffM\times\DiffPi \ | \ \difP\circ\mrsfunc=\mrsfunc\circ\difM \}$ be the corresponding stabilizer and  
$\OrbfMPi=\{\difP\circ\mrsfunc\circ\difM^{-1} \ | \ (\difM,\difP)\in\DiffM\times\DiffPi\}$ the orbit of $\mrsfunc$.

Evidently, 
$\StabfM \equiv \StabfM\times\id_{\Psp}\subset\StabfMPi$ and $\OrbfM\subset\OrbfMPi.$

The aim of this paper is to show that for almost all mappings $\mrsfunc\in\smone$ and some natural subgroups $\DiffPi\subset\DiffP$ we have the following homotopy equivalences (in the corresponding $C^{\infty}$-topologies):
$\StabfM\approx\StabfMPi$, 
$\OrbfM\approx\OrbfMRi$ for $\Psp=\RRR$, and 
$\OrbfM\times\Circle\approx\OrbfMSi$ for $\Psp=\Circle$.
In fact, we obtain precise relationships between the topological (not just homotopy) types of these spaces.

\subsection{Conditions on $\mrsfunc$.}
In order to formulate our main results (Theorems~\ref{th:Stab} and~\ref{th:Orbit}) we introduce the following conditions \condBndCrVal\ and \condJacId.
Say that $\mrsfunc\in\smone$ satisfies the condition \condBndCrVal\ if
\begin{enumerate}
\item[\condBndCrVal]
$\mrsfunc$ is constant at every  connected component of $\partial\manif$ and has only finitely many \emph{critical values}.
\end{enumerate}

For each point $\pnt\in\manif$ let $\germm$ be the algebra of germs of smooth functions at $\pnt$.
If $\mrsfunc\in\germm$, then the {\em Jacobi ideal\/} $\Jacidfz$ of $\mrsfunc$ at $\pnt$ is the ideal in $\germm$ generated by partial derivatives of $\mrsfunc$ at $\pnt$.
Evidently, it does not depend on a particular choice of local coordinates at $\pnt$.

We will say that a mapping $\mrsfunc\in\smone$ satisfies the condition \condJacId\ if the following holds true.
\begin{enumerate}
\item[\condJacId]
Let $\pnt$ be a critical point of $\mrsfunc$ and $\mrsfunc:\Rm\to\RRR$ a local representation of $\mrsfunc$ at $\pnt$ such that $\mrsfunc(\pnt)=0$.  
Then the germ of a function $\mrsfunc$ at $\pnt$ belongs to the Jacobi ideal $\Jacidfz$.
\end{enumerate}

This means that there are smooth functions $\Fld_{1},\ldots,\Fld_{\dimM}$ such that
\begin{equation}\label{equ:f_in_Jacid}
\mrsfunc(x) =
\sum_{i=1}^{\dimM} \mrsfunc'_{x_{i}}(x) \Fld_{i}(x), \qquad x=(x_1,\ldots,x_{\dimM})\in\Rm.
\end{equation}

Moreover, if we define a vector field $\Fld$ near $\pnt$ by
$\Fld =(\Fld_{1},\ldots,\Fld_{\dimM})$, then~\eqref{equ:f_in_Jacid} can also be written in the following form:
\begin{equation}\label{equ:f_in_Jacida_verct_fld}
\mrsfunc(x)  = \dfvalx{\mrsfunc}{\Fld}{x}, \qquad x\in\RRR^{\dimM}
\end{equation}
where $\dfval{\mrsfunc}{\Fld}$ is a derivative of $\mrsfunc$ along $\Fld$.

Notice, that conditions \condBndCrVal\ and \condJacId\ are invariant with respect to left-right action.

\subsubsection{}
Suppose that a function $\mrsfunc\in\smone$ satisfies \condBndCrVal.
Then the set of {\em critical points\/} of $\mrsfunc$ may be infinite.
Moreover, there may be critical points on $\partial\manif$.

The values of $\mrsfunc$ on the connected components of $\partial\manif$ will be called \emph{boundary\/} ones.
All critical and boundary values of $\mrsfunc$ will be called \emph{exceptional}
and the inverse images of these exceptional values under $\mrsfunc$ will be called \emph{exceptional} level-sets.
Since $\manif$ is assumed compact, it follows that the set of exceptional values is finite.

Let $\nv$ be the total number of exceptional values of $\mrsfunc$.

If $\nv=0$, then $\manif$ is closed, $\Psp=\Circle$, and $\mrsfunc:\manif\to\Circle$ is a locally trivial fibration.

Otherwise, $\nv\geq 1$.
If $\Psp=\Circle$, then we shall always regard $\Circle$ as the group of reals modulo $\nv$:
$$ \Circle \ \equiv \ \RRR/\nv\ZZZ, $$ 
not $\RRR/\ZZZ$ as usual!
Therefore in both cases of $\Psp$ we can assume that $1,\ldots,\nv$ are all of the exceptional values of $\mrsfunc$.
But for $\Psp=\Circle$, they are taken modulo $\nv$, and in particular we have that $\nv\equiv0$. 

\smallskip

Let us define the following groups.
If $\mrsfunc\in\smr$, then let

$\bullet$~$\DiffRf$ be the subgroup of $\DiffR$ consisting of diffeomorphisms that preserve orientation of $\RRR$, have compact support, and leave the image $[1,\nv]$ of $\mrsfunc$ invariant;

$\bullet$~$\DiffRfe$ be the subgroup of $\DiffRf$ consisting of diffeomorphisms that also fix every exceptional value $1,\ldots,\nv$ of $\mrsfunc$;

$\bullet$~$\DiffMRf = \DiffM \times \DiffRf$.

\smallskip

If $\mrsfunc\in\sms$, then let

$\bullet$~$\DiffSpl$ be the group of preserving orientation diffeomorphisms of $\Circle$;

$\bullet$~$\DiffSfE$ be the subgroup of $\DiffSpl$ preserving the set $\{1,\ldots\nv\}$ of exceptional values of $\mrsfunc$. 
If $\nv=0$, then $\DiffSfE=\DiffSpl$;

$\bullet$~$\DiffSfe$ be the (normal) subgroup of $\DiffSfE$ fixing the set $\{1,\ldots\nv\}$ point-wise, thus
$\DiffSfE/\DiffSfe$ is a cyclic group $\ZZZ_{\nv}$ of order $\nv$;

$\bullet$~$\DiffMSf = \DiffM \times \DiffSpl$.

\smallskip

Then $\DiffM$ and $\DiffMPf$ act on $\smone$ by formulas~\eqref{equ:act_M}
and~\eqref{equ:act_MP}. Let $\StabfM$, $\StabfMP$, $\OrbfM$, and $\OrbfMP$ be respectively the stabilizers and the orbits of $\mrsfunc$ under these actions.
Evidently,
$$ \StabfM\times\id_{\Psp}\subset\StabfMP \qquad\text{and}\qquad
\OrbfM\subset\OrbfMP. $$
Finally, we endow the spaces $\DiffM$, $\DiffP$, and $\smone$ with the corresponding $C^{\infty}$ Whitney topologies.
These topologies yield certain topologies on $\DiffMPf$ and on the corresponding stabilizers and orbits of $\mrsfunc$.

\subsection{Stabilizers.}
Let us say that a diffeomorphism $\difP\in\DiffP$ is {\em left-trivial\/} or {\em \LT-trivial\/} for $\mrsfunc$ if there exists a diffeomorphism $\difM\in\DiffM$ such that $(\difM,\difP)\in\StabfMP$, i.e. $\difP\circ\mrsfunc=\mrsfunc\circ\difM$.
Thus applying $\difP$ to $\mrsfunc$ (acting from the left) we remain in the {\em right\/} orbit $\OrbfM$ of $\mrsfunc$.
This explains the term ``left-trivial''.

Let $\prMPP:\DiffM\times\DiffP\to\DiffP$ be the standard projection. 
Evidently, $\prMPP$ is a homomorphism.
Consider the restriction of $\prMPP$ to $\StabfMP$, then its kernel coincides with $\StabfM$, and the image $\prMPP(\StabfMP) \subset\DiffP$ consists of all \LT-trivial for $\mrsfunc$ diffeomorphisms.

\begin{thm}\label{th:Stab}
Suppose that $\mrsfunc\in\smone$ satisfies \condBndCrVal\ and \condJacId.
Then  
$$
\begin{array}{lcl}
\text{{\rm(1)}~Case $\Psp=\RRR$:} & \qquad & \prMPP(\StabfMR)\,=\,\DiffRfe, \\ [2.5mm]
\text{{\rm(2)}~Case $\Psp=\Circle$, $\nv\geq1$:} & \qquad & \DiffSfe \, \subseteq \, \prMPP(\StabfMS)\, \subseteq \, \DiffSfE. \\[2.5mm]
\text{{\rm(3)}~Case $\Psp=\Circle$, $\nv=0$:} & \qquad & \prMPP(\StabfMS)\,=\,\DiffSpl.
\end{array}
$$
In the cases (1) and (2) $\prMPP$ admits a continuous section $\sectMPP:\DiffPfe\to\StabfMP$ i.e. $\prMPP\circ\sectMPP=\id(\DiffPfe)$. Moreover, this section is a homomorphism.

In the case (3) $\prMPP$ admits a continuous section $\sectMPP:\DiffSpl\to\StabfMS$ iff the fibration $\mrsfunc:\manif\to\Circle$ trivial. Such a section can also be chosen to be a homomorphism.
\end{thm}
\begin{rem}\label{rem:projection}
The section of the projection $\prMPP$ in Theorem~\ref{th:Stab} must be of the form $\sectMPP(\difR)=(\difShift(\difP),\difP)$, 
where $\difShift:\DiffPfe\to\DiffM$ is a continuous mapping such that for every diffeomorphism $\difP\in\DiffPfe$ (fixing each exceptional value of $\mrsfunc$) we have
$\difP\circ\mrsfunc = \mrsfunc \circ \difShift(\difP)\in\OrbfM$. 
Notice also that $\sectMPP$ is a homomorphism iff $\difShift$ is.
\end{rem}

\begin{rem}
In the case (2) 
denote $\pStabfMS=\prMPP^{-1}(\DiffSe)$.
Since $\DiffSfE/\DiffSfe\approx\ZZZ_{\nv}$, it follows from Theorem~\ref{th:Stab}
that 
\begin{equation}\label{equ:pStabfMS_DiffSfe_Zc}
\StabfMS\,/\,\pStabfMS \ \approx \
\prMPP(\StabfMS)\,/\,\DiffSfe\  \approx\ \ZZZ_{\kcnv},
\end{equation}
for some $\kcnv$ that divides $\nv$.
\end{rem}

\begin{ThStabA}[Another formulation]
The following sequences of group homomorphisms are exact
$$
\begin{array}{lcl}
\text{{\rm(1)}~Case $\Psp=\RRR$:} & \qquad &
1 \ \to \ \StabfM \ \to \ \StabfMR \ \xrightarrow{\prMPP} \ \DiffRfe \ \to \ 1,
 \\ [2mm]
\text{{\rm(2)}~Case $\Psp=\Circle$, $\nv\geq1$:} & \qquad &
1 \ \to \ \StabfM \ \to \ \pStabfMS \ \xrightarrow{\prMPP} \ \DiffSfe \ \to \ 1,
\\ [2mm]
\text{{\rm(3)}~Case $\Psp=\Circle$, $\nv=0$:} & \qquad &
1 \ \to \ \StabfM \ \to \ \StabfMS \ \xrightarrow{\prMPP} \ \DiffSpl \ \to \ 1.
\end{array}
$$
They always split in the cases (1) and (2), and in the case (3) iff $\mrsfunc:\manif\to\Circle$ is a trivial fibration.
\end{ThStabA}

Notice that the existence of splittings in the cases (1) and (2) is rather natural, since $\prMPP$ is a principal $\StabfM$-fibration and $\DiffRfe$ and $\DiffSfe$ for $n\geq1$ are contractible, see Lemma~\ref{lm:DiffP_is_contr}.

Thus in the case (1) we obtain a homeomorphism
$\StabfMR \cong \StabfM\times\DiffRfe$, whence the embedding $\StabfM\subset\StabfMR$ is a homotopy equivalence.

In the case (2) we have that $\pStabfMS\cong\StabfM\times\DiffSfe$.
Since $\StabfMS/\pStabfMS$ is a cyclic group $\ZZZ_{\kcnv}$, we get 
$\StabfMS\cong\StabfM\times\DiffSfe\times\ZZZ_{\kcnv}$.
Whence $\StabfMS$ is homotopy equivalent to $\StabfM\times\ZZZ_{\kcnv}$.

Moreover, suppose that $\ZZZ_{\kcnv}$ is non-trivial, i.e. $\kcnv>1$. 
Then there exists $(\difM,\difP)\in\StabfMS$, i.e. $\difP\circ\mrsfunc=\mrsfunc\circ\difM$, such that $\difP$ cyclically shifts exceptional values $1,\ldots,\nv$ of $\mrsfunc$ and $\difM$ cyclically shifts the corresponding exceptional level-sets of $\mrsfunc$:
$$
\difP(\indx) = \indx+n/\kcnv, 
\qquad 
\difM(\exlevel_{\indx}) = \exlevel_{\indx+n/\kcnv}, \qquad (\indx=1,\ldots,n),
$$
where $\exlevel_{\indx}=\mrsfunc^{-1}(\indx)$ and each sum is taken modulo $n$.
Thus the level-sets 
$$\exlevel_{\indx}, \ \exlevel_{\indx+n/\kcnv}, \ \ldots \ \exlevel_{\indx+n(\kcnv-1)/\kcnv}$$
are homeomorphic each with other. 
This situation is not typical, though it could be stable (e.g.\!\! for generic Morse function).
Thus for most functions we should have $\DiffSfe\,=\,\prMPP(\StabfMS)$.
In this case $\StabfMS$ is homeomorphic with $\StabfM\times\DiffSfe$ and 
the embedding $\StabfM\subset\StabfMS$ is a homotopy equivalence.

\subsubsection{Interpretation: holonomy.}
Let $\mrsfunc:\manif \to B$ be a finite-dimensional vector bundle over a smooth manifold $B$.
For each $b\in B$ let $\manif_{b}=\mrsfunc^{-1}(b)$ be the corresponding fiber.
Choose some connection on $\manif$.
Then for each smooth path $\omega:I\to B$ there exists a smooth isotopy (consisting even of linear isomorphisms) 
$\difM_t:\manif_{\omega(0)}\to \manif_{\omega(t)} \subset \manif,$
called {\em holonomy} along $\omega$.
It follows that, if $\difP_t:B \to B$ is an isotopy and $\difP_{0}=\id_{B}$, then there is an isotopy $\difM_t:\manif\to\manif$ 
such that $\difM_0=\id_{\manif}$ and the following diagram is commutative:
$$
\begin{CD}
\manif @>{\mrsfunc}>> B \\
@V{\difM_t}VV @VV{\difP_t}V \\
\manif @>{\mrsfunc}>> B
\end{CD}
\qquad \text{i.e.}  \qquad \difP_t\circ\mrsfunc=\mrsfunc\circ\difM_t.
$$
In other words, $(\difM_t, \difP_t)$ belongs to the stabilizer of $\mrsfunc$ under the left-right action of the group $\DiffM\times\Diff_{B}$ on $C^{\infty}(\manif,B)$.
Moreover, if $G\subset\Diff_{B}$ is a {\em simply-connected} subgroup, then we obtain a continuous homomorphism $\difShift:G \to \DiffM$ such that $\difP\circ\mrsfunc=\mrsfunc\circ\difShift(\difP)$ for all $\difP\in G$.
   
Notice that the vector-bundle projection $\mrsfunc$ is a smooth mapping {\em without\/} critical points.
If we replace $\mrsfunc$ with an arbitrary smooth mapping between arbitrary smooth manifolds, and try to construct such a "holonomy homomorphism", then we should put some restrictions on $\mrsfunc$ and $G$. In particular, $G$ must preserve the image of $\mrsfunc$ and its ``exceptional'' values.
Theorem~\ref{th:Stab} gives rather general conditions when such a homomorphism exists for the case $\dim B=1$. 

\subsection{Orbits.}
We will now describe the relationships between the orbits.
\begin{defn}\label{defn:essent}
A critical point $\pnt$ of $\mrsfunc\in\smone$ is \emph{essential}, if for
every neighborhood $\nbh$ of $\pnt$ there is a neighborhood
$\Nbh$ of $\mrsfunc$ in $\smone$ with $C^{\infty}$-topology such
that each $\gfunc\in\Nbh$ has a critical point in $\nbh$.
\end{defn}
\begin{exmp}
Let $\mrsfunc(x)=x^2$ and $\gfunc(x)=x^3$. 
Then $0\in\RRR$ is an essential critical point for $\mrsfunc$ but not for $\gfunc$.
\end{exmp}

\begin{thm}\label{th:Orbit}
Suppose that $\mrsfunc$ satisfies the conditions \condBndCrVal\ and \condJacId\ and 
each {\em critical\/} level-set of $\mrsfunc$ includes either an essential critical point or a connected component of $\partial\manif$. 

If $\Psp=\RRR$, then the embedding $\OrbfM\subset\OrbfMR$ extends to a homeomorphism
$$\OrbfM \ \times \ \RRR^{n-2} \approx \OrbfMR.$$

Let $\Psp=\Circle$ and $\kcnv$ be the index of $\DiffSfe$ in $\prMPP(\StabfMS)$, as in Theorem~\ref{th:Stab}.
\begin{enumerate}
\item[{\rm a)}] If $\nv=0$, then $\OrbfM\,=\,\OrbfMS$.

\item[{\rm b)}] If \ $\nv$ \ is even and \ $\nv/\kcnv$ \ is odd, then 
$$\OrbfM \,\times\,  \Circle \,\tilde\times\, \RRR^{\nv-1} \,\approx\, \OrbfMS,$$
 where $\Circle\tilde\times\RRR^{\nv-1}$ is the total space of a unique non-trivial $(\nv-1)$-dimensional fibration over $\Circle$.
 
\item[{\rm c)}] Otherwise, 
	$$\OrbfM  \,\times\,  \Circle  \,\times\, \RRR^{\nv-1} \,\approx\, \OrbfMS.$$
\end{enumerate}
\end{thm}

\subsection{Structure of the paper.}
In Section~\ref{sect:grp_Vdfa} for each germs of ``admissible'' smooth functions $\afunc:\RRR\to\RRR$ at $0\in\RRR$ (see Definition~\ref{def:admiss_func}) we introduce and study a certain group $\Vdf{\afunc}$ of germs of diffeomorphisms $\RRR\to\RRR$ at $0$. 

In Section~\ref{sect:local_act} we consider the local left-right action of the groups of germs of diffeomorphisms of $\Rm$ and $\RRR$ on the germs of smooth functions $\mrsfunc:\Rm\to\RRR$ such that $\mrsfunc(0)=0$. 
We give a sufficient condition \condJacIda{\afunc}\ on $\mrsfunc$ when $\Vdf{\afunc}$ consists of \LT-trivial for $\mrsfunc$ diffeomorphisms 
 (Theorem~\ref{th:ex_loc_sect}).
The most complete result (Theorem~\ref{th:loc_Stab}) which is also a local variant of Theorem~\ref{th:Stab} is obtained for $\mrsfunc$ satisfying condition \condJacIda{\id}\ being a local analogue of \condJacId.

Section~\ref{sect:cond_JacId}. 
We show that the condition \condJacId\ holds for a very large class of singularities and is invariant with respect to a stable equivalence of singularities (Lemma~\ref{lm:stable_equiv}).
On the other hand, there are singularities that do not satisfy this condition (Claim~\ref{clm:ex_not_in_J}).

In Section~\ref{sect:proofA} we prove Theorem~\ref{th:Stab}.

Section~\ref{sect:DiffPfe}. The finite-dimensional spaces of adjacent classes of $\DiffRf/\DiffRfe$ and $\DiffSf/\DiffSfe$ (Theorem~\ref{th:struct_DP_DPe}) are described.

Section~\ref{sect:cont_crvmap}.
We give a sufficient condition for the continuity of the mapping $\crvmap$ corresponding to each $\gfunc\in\OrbfMP$ the ordered set of exceptional values of $\gfunc$ (Lemma~\ref{lm:cond_crlev_cont}).
We also show that without this condition $\crvmap$ may loose continuity.

Finally in Section~\ref{sect:proofB} we prove Theorem~\ref{th:Orbit}.

\section{Groups $\Vdf{\afunc}$}\label{sect:grp_Vdfa}
Let $\smzr$ be the algebra of germs of smooth functions at $0\in\RRR$.
For each $\mfunc\in\smzr$ we will denote by $\mideal{\mfunc}$ the ideal $\mfunc\cdot\smzr$ in $\smzr$. 

\begin{defn}\label{def:admiss_func}
We will say that a function $\afunc\in\smzr$ is {\em admissible}, provided
\begin{enumerate}
\item
$\afunc(0)=0$, $\afunc'(0)$ is equal either $0$ or $1$, and 
\item
there is a neighborhood $\nbh$ of $0\in\RRR$ such that the intersection $\afunc^{-1}(0)\cap\nbh$ is {\em nowhere dense in $\nbh$}, thus $\afunc$ in not constant on open intervals near $0$. 
\end{enumerate}
\end{defn}
For each admissible $\afunc\in\smzr$ we will now define a certain group $\Vdf{\afunc}$ of germs of diffeomorphisms of $\RRR$ at $0$.
They are analogous to the groups $G_{d}$ of~\cite[\S5]{Arnold74} of quasi-homogeneous diffeomorphisms $\CCC^n\to\CCC^n$ of order $d\geq0$. 
Our situation is simpler since we consider diffeomorphisms of $\RRR$, on the other hand dimension $1$ allows to prove more.
The proximity to the identity $\id_{\RRR}$ for diffeomorphisms of $\Vdf{\afunc}$ is defined not just {\em up to order $d$\/} but {\em up to the admissible function $\afunc$\/} which can be flat.
Moreover, our approach to $\Vdf{\afunc}$ differs from~\cite[\S5]{Arnold74}.
We study this group using its characterization as the set of smooth shifts along trajectories of the vector field $\afunc(\vars)\frac{d}{d\vars}$ on $\RRR$, see Theorem~\ref{th:Vdf_shifts}.
The formulas for shift-functions will play a key role in the proof of Theorem~\ref{th:Stab} and its local variant Theorem~\ref{th:ex_loc_sect}.

\subsection{Definition of $\Vdf{\afunc}$.}\label{sect:def_Vdfa}
Let $\afunc\in\smzr$ be an admissible function. 
Then we define $\Vdf{\afunc}$ to be the the subset of $\DiffRz$ consisting of {\em preserving orientation diffeomorphisms\/} $\adif$ of the following form:
\begin{equation}\label{equ:reps_difR_all}
\adif(s)= s + \afunc(s) \adgfn(s), \qquad \adgfn\in\smzr.
\end{equation}
In other words, \  $\adif - \id_{\RRR}  \in \mideal{\afunc}$.
Consider two cases.

\caseA~Suppose that $\afunc(0)=0$ and $\afunc'(0)=1$, i.e. $\afunc(\vars)=\vars\bafunc(\vars)$, where $\bafunc\in\smzr$ and $\bafunc(0)=\afunc'(0)=1$.
In this case $\Vdf{\afunc}=\DiffRz$.
Indeed, if $\difR\in\DiffRz$, then $\difR(0)=0$, whence 
$\difR(\vars)-\vars = \vars\ofunc(\vars) = \afunc(\vars)\frac{\ofunc(\vars)}{\bafunc(\vars)}$, where $\ofunc\in\smzr$.
Thus $\difR\in\Vdf{\afunc}$.

\caseB~Otherwise, $\afunc(0)=\afunc'(0)=0$, i.e.\!\! $\afunc(\vars)=\vars^2\bafunc(\vars)$ for some $\bafunc\in\smzr$.
Then each $\adif\in\Vdf{\afunc}$ has the following form $\adif(\vars)=\vars+\vars^2\bafunc(\vars)\adgfn(\vars)$.
In particular, $\adif'(0)=1$.

Notice that $\Vdf{\vars^{k}}$ consists of diffeomorphisms of the form 
$\adif(\vars)=\vars+\vars^{k}\adgfn(\vars)$.
Equivalently, $\adif\in\Vdf{\vars^{k}}$ iff $\adif'(0)=1$, and $\adif^{(p)}(0)=0$ for $p=2,3,\ldots,k-1$.

\begin{rem}
$\Vdf{\afunc\cfunc}\subseteq\Vdf{\afunc}$ for each $\cfunc\in\smzr$,
and $\Vdf{\afunc\cfunc}=\Vdf{\afunc}$ iff $\cfunc(0)\not=0$.
In particular, $\Vdf{\afunc}\subset\Vdf{\id}$ for all admissible $\afunc$.
\end{rem}
\begin{rem}\label{rem:adgfn_is_det_by_adif}
Condition (2) of Definition~\ref{def:admiss_func} implies that $\adgfn$ is uniquely determined near $0$ by $\adif$. 
Indeed, if this condition is violated, then there is a sequence of mutually disjoint closed intervals $A_k$ converging to $0$ and such that $\afunc|_{A_k}\equiv0$.
Then varying $\adgfn$ on $A_k$ we do not change $\adif$.
\end{rem} 	

\subsection{Another description of $\Vdf{\afunc}$}\label{sect:loc_flow}
We show that $\Vdf{\afunc}$ coincides with the set of smooth shifts along trajectories of the vector field $\afunc(\vars)\frac{d}{ds}$ on $\RRR$. 

Let $\afunc:\rnbh\to\RRR$ be an admissible smooth function. 
Define a vector field $\rFld$ on $\RRR$ by $\rFld(\vars)=\afunc(\vars)\frac{d}{d\vars}$.
Let also $\rflow:\rnbh\times\Interv\to\RRR$ be the local flow generated by $\rFld$, where $\rnbh$ is a neighborhood of $0\in\RRR$ and $\Interv$ an open interval containing $0\in\RRR$.

\begin{thm}\label{th:Vdf_shifts}
Let $\adif\in\smzr$.
Then $\adif\in\Vdf{\afunc}$  iff $\adif$ is a smooth shift along trajectories of $\rflow$, i.e.
$\adif(\vars)=\rflow(\vars,\ashfn(\vars))$ for some $\ashfn\in\smzr$.

Moreover $\adgfn=\ashfn\,\ofunc$, where $\ofunc\in\smzr$ and $\ofunc(0)=1$. 
Hence $\adgfn(0)=\ashfn(0)$.
\end{thm}
We will call $\ashfn$ a {\em shift-function\/} of $\adif$ with respect to $\rflow$.
Before proving this theorem let us deduce some corollaries.

\begin{lem}\label{lm:Vdf_is_grp}
{\rm (c.f.~\cite[Prop.~5.2]{Arnold74})}.
$\Vdf{\afunc}$ is a group.
\end{lem}
\proof
Let $\adif,\bdif\in\Vdf{\afunc}$.
By Theorem~\ref{th:Vdf_shifts} $\adif(\vars)=\rflow(\vars,\ashfn(\vars))$ and $\bdif(\vars)=\rflow(\vars,\bshfn(\vars))$ for some $\ashfn,\bshfn\in\smzr$.
It is easy to see \cite[Proposition~3]{Maks:Shifts} that 
\begin{equation}\label{equ:shift_func_formulas}
\bdif\circ\adif(\vars)=\rflow\bigl(\vars, \ashfn(\vars) + \bshfn\circ\adif(\vars) \bigl), \qquad 
\bdif^{-1}(\vars) = \rflow\bigl(\vars, - \bshfn\circ\bdif^{-1}(\vars)\bigr).
\end{equation}
Hence, again by Theorem~\ref{th:Vdf_shifts} \ \ $\bdif\circ\adif, \ \bdif^{-1} \in \Vdf{\afunc}$. 
Thus $\Vdf{\afunc}$ is a group.
\endproof

\begin{lem}\label{lm:adif_in_Vdf_am}
Let \ $\adif\in\Vdf{\afunc}$ \ and \ $\mfunc\in\smzr$ \ be such that \ $\afunc\mfunc\in\smzr$ \ is admissible. 
Then the following conditions are equivalent:
$$ 
{\rm (1)}~\adif\in\Vdf{\afunc\mfunc} \qquad\qquad
{\rm (2)}~\adgfn\in\mideal{\mfunc} \qquad\qquad
{\rm (3)}~\ashfn\in\mideal{\mfunc}.
$$
Thus the group $\Vdf{\afunc\mfunc}$ consists of all smooth shifts along trajectories of $\rflow$ whose shift-functions are proportional to $\mfunc$, i.e. belong to the ideal $\mideal{\mfunc}$.
\end{lem}
\proof
The equivalence (2)$\Leftrightarrow$(3) holds by Theorem~\ref{th:Vdf_shifts}, since $\mideal{\adgfn}=\mideal{\ashfn}$.

(1)$\Leftrightarrow$(2)~
Notice that $\adif\in\Vdf{\afunc\mfunc}$ means that $\difR(\vars)=\vars+\afunc(\vars)\mfunc(\vars)\ofunc(\vars)$ for some $\ofunc\in\smzr$.
By~\eqref{equ:reps_difR_all} and Remark~\ref{rem:adgfn_is_det_by_adif} this condition is equivalent to the following
\ $\adgfn\,=\,\mfunc\,\ofunc\,\in\,\mideal{\mfunc}$.
\endproof

\begin{lem}
Let $\bdif,\adif\in\Vdf{\afunc}$ and $\cdif=\bdif\circ\adif\circ\bdif^{-1}$.
Then $ \cshfn = \ashfn \circ \bdif^{-1} \cdot \nfunc,$ \ where $\nfunc\in\smzr$ and $\nfunc(0)=1$.
Hence \ $\cdif\in\Vdf{\afunc\mfunc}$ \ iff \ $\ashfn\circ\bdif^{-1} \in \mideal{\mfunc}$.
\end{lem}
\proof 
From~\eqref{equ:shift_func_formulas} we get  
$$
\cdif(\vars) = \bdif\circ\adif\circ\bdif^{-1}(\vars)=
\rflow\bigl(\vars, 
								-\bshfn\circ\bdif^{-1}(\vars) + 
								\ashfn\circ\bdif^{-1}(\vars) + 
								\bshfn\circ\adif\circ\bdif^{-1}(\vars)
			\bigl).
$$
Thus \ $\cshfn=(\ashfn + \bshfn\circ\adif - \bshfn)\circ\bdif^{-1}$.

Then it follows from Hadamard lemma and Theorem~\ref{th:Vdf_shifts} that 
$$
\bshfn\circ\adif(\vars) - \bshfn(\vars) \, = \,
(\adif(\vars)-\vars)\,\,\bar\bshfn\, = \,
\afunc\,\adgfn\,\bar\bshfn \, 
= \, 
\afunc\,\ashfn\,\ofunc\,\bar\bshfn
$$ 
for some $\bar\bshfn\in\smzr$.
Hence \
$\cshfn = 
\bigl( \ashfn \cdot (1\,+\, \afunc\,\ofunc\, \bar\bshfn ) \bigr) \circ \bdif^{-1} = 
\ashfn\circ\bdif^{-1}\cdot\nfunc,$ \ 
where 
$\nfunc=(\id_{\RRR}\, +\, \afunc\,\ofunc\, \bar\bshfn ) \circ \bdif^{-1}$ and $\nfunc(0)=1$.
\endproof

\begin{cor}\label{cor:Vdf_am_normal_char}
$\Vdf{\afunc\mfunc}$ is a normal subgroup of $\Vdf{\afunc}$ iff the ideal $\mideal{\mfunc}\subset\smzr$ is invariant with respect to the action of $\Vdf{\afunc}$ on $\smzr$ by \
$$\bdif\cdot\sfunc = \sfunc\circ\bdif^{-1}, \qquad \bdif\in\Vdf{\afunc}, \ \sfunc\in\smzr.$$
\end{cor}
\begin{cor}\label{cor:Vdf_ska_is_normal_in_Vdf_a}
{\rm (c.f.~\cite[Prop.~5.3]{Arnold74})}.
For each $k\geq1$ the group $\Vdf{\vars^{k}\afunc}$ is normal in $\Vdf{\afunc}$.
\end{cor}
\proof
It suffices to show that $\sfunc\circ\adif\in\mideal{\vars^{k}}$
for every $\sfunc\in\mideal{\vars^{k}}$ and $\adif\in\Vdf{\id}$.
Indeed, we have that $\sfunc(\vars)=\vars^{k}\bar\sfunc$ and 
$\adif(\vars) = \vars\,\ofunc(\vars)$ for some $\bar\sfunc,\ofunc\in\smzr$.
Then
$$
\sfunc\circ\adif(\vars) \,=\, \adif(\vars)^{k}\,\bar\sfunc(\adif(\vars)) \,=\, 
 \vars^{k}\,\ofunc(\vars)^{k}\,\bar\sfunc(\adif(\vars))\,\in\mideal{\vars^{k}}.
\qed
$$

Define the following mapping
$\hmV{\afunc}:\Vdf{\afunc}\to\RRR$ by 
$\hmV{\afunc}(\adif) = \adgfn(0)=\ashfn(0)$.
Then it follows from~\eqref{equ:shift_func_formulas} that $\hmV{\afunc}$ is a surjective homomorphism whose kernel is $\Vdf{\vars\afunc}$.

Thus for every admissible $\afunc$ we obtain a sequence of normal subgroups 
$$
\Vdf{\afunc} \ \supset \ \Vdf{\vars\afunc} \ \supset \ \Vdf{\vars^2\afunc} \ \supset \ \cdots
$$
such that each factor is isomorphic to $\RRR$.
Hence for each $k$ the group $\Vdf{\afunc}/\Vdf{\vars^{k}\afunc}$ is a Lie group diffeomorphic with $\RRR^{k}$.
{\rm (c.f.~\cite[Prop.~5.5]{Arnold74})}.

\subsection{Proof of Theorem~\ref{th:Vdf_shifts}.}
It suffices to establish the following proposition that describes formulas for $\rflow$.
\begin{prop}\label{pr:char_flowR1}
There is a smooth function $\cfunc$ on $\rnbh\times\Interv$ such that 
\begin{equation}\label{equ:rflow_2}
\rflow(\vars,t) = \vars + t \afunc(\vars)\cfunc(\vars,t),
\end{equation}
and $\cfunc(0,t)\equiv1$.
In particular, $\rflow_t \in \Vdf{\afunc}$ for every $t\in\Interv$.
\end{prop}

First we deduce Theorem~\ref{th:Vdf_shifts}.

{\bf Sufficiency.}
It follows from~\eqref{equ:rflow_2} that if 
$$\adif(\vars)\,=\,\rflow(\vars,\ashfn(\vars))\,=\,\vars\,+\,\afunc(\vars)\,\ashfn(\vars)\,\cfunc(\vars,\ashfn(\vars)),$$
then $\adif\in\Vdf{\afunc}$ with \ $\adgfn(\vars)=\ashfn(\vars)\,\cfunc(\vars,\ashfn(\vars))$ \ and \
$\ofunc(\vars)=\cfunc(\vars,\ashfn(\vars))$, where 
$\ofunc(0)=\cfunc(0,\ashfn(0))=1$.

{\bf Necessity.}
Suppose that $\adif\in\Vdf{\afunc}$.
In order to show that $\adif(\vars)=\rflow(\vars,\ashfn(\vars))$ we have to show that 
$\adgfn(\vars)=\ashfn(\vars)\cfunc(\vars,\ashfn(\vars))$ for some $\ashfn\in\smzr$.

Consider the function $\dfunc(\vars,t)=t\cfunc(\vars,t)$.
Notice that $\dfunc'_{t}(\vars,0)=\cfunc(\vars,0)=1$, whence there exists a smooth function $\qfunc(\vars,t)$ such that $t = \dfunc(\vars,\qfunc(\vars,t))$.
Therefore, we put $\ashfn(\vars) = \qfunc(\vars,\adgfn(\vars))$. 
Then $\adgfn(\vars)=\dfunc(\vars,\ashfn(\vars))=\ashfn(\vars)\cdot\cfunc(\vars,\ashfn(\vars))$. 
\endproof

\subsection{Proof of Proposition~\ref{pr:char_flowR1}.} 
For simplicity we will sometimes omit the dependence on $\vars$ and $(\vars,t)$.
Recall that $\rflow'_{t}(\vars,0)=\afunc(\vars)$, whence the Taylor expansion of $\rflow(\vars,t)$ in $t$ at $(\vars,0)$ has the following form:
\begin{equation}\label{equ:flow_decomp}
\rflow(\vars,t) \,= \,\rflow(\vars,0) \,+\,  t \,\rflow'_t(\vars,0)\,
 + \,\cdots\, = \, \vars \,+ \,t\,\afunc(\vars)\, + \,\cdots
\end{equation}
Then 
$$
\cfunc(\vars,t) = \frac{\rflow(\vars,t)-\vars}{t\,\afunc(\vars)}.
$$
This function is defined only for those $(\vars,t)$ for which $t\afunc(\vars)\not=0$.
Nevertheless, since $\rflow(s,0)=s$, it follows that $(\rflow(s,t)-s)/t$ is smooth.
Moreover, for $t\not=0$ we have that $\afunc(\vars)=0$ iff $\rflow(\vars,t)=\vars$.
\begin{lem}\label{lm:cfunc_rel}
$\cfunc(\vars,t)$ satisfies the following differential equation:
\begin{equation}\label{equ:dif_eq_cfunc}
\cfunc\,'_{\vars}(\vars,t) = \cfunc^2(\vars,t) \cdot t \cdot \mu(\vars,t),
\end{equation}
where $\mu$ is a certain smooth function on $\nbh\times\Interv$.
Whence 
$$\cfunc(\vars,t)=\frac{1}{c(t) \, - \, t \, \int\limits_{0}^{\vars} \mu(z,t) dz},$$
where $c(t)$ is a smooth function such that $c(0)=1$.
Thus $\cfunc$ is smooth on $\nbh\times(-\eps,\eps)$ for sufficiently small $\eps>0$ and $\cfunc(0,t)=\frac{1}{c(0)}=1$.
\end{lem}
\proof 
A simple calculation shows that if $\afunc(\vars)\not=0$ then 
\begin{equation}\label{equ:cfunc_main_rel}
\cfunc\,'_{\vars}(\vars,t) = \frac{\afunc\cdot\rflow'_{\vars}\, - \, \afunc \, - \, (\rflow-\vars) \, \afunc\,'_{\vars}}{t\,\afunc^2}.
\end{equation}

\begin{claim}
The first term of the numerator in~\eqref{equ:cfunc_main_rel} is equal to 
$$\afunc(\vars) \,\cdot \,\rflow'_{\vars}(\vars,t) \, = \, \afunc\circ\rflow(\vars,t).$$
\end{claim}
\proof
Notice that $\rFld$ defines the following differential equation on $\RRR$:
$\frac{d\vars}{dt} = \afunc(\vars)$, whence $dt = \frac{d\vars}{\afunc(\vars)}$.
Then for every $\vars\in\rnbh_1$ the time $t$ along the trajectory of $\rflow$ between $\vars$ and $\rflow(\vars,t)$ is equal to \
$ t = \int\limits_{0}^{t} dt =  \int\limits_{\vars}^{\rflow(\vars,t)} \frac{dz}{\afunc(z)}$,
(notice that if $\afunc(\vars)\not=0$, then $\afunc\not=0$ between $\vars$ and $\rflow(\vars,t)$\,). 

Differentiating both sides of this equality in $\vars$ we get 
$$
0 = \frac{\rflow'_{\vars}}{\afunc(\rflow)} - \frac{1}{\afunc} =
\frac{\afunc  \cdot \rflow'_{\vars} - \afunc(\rflow)}{\afunc \cdot \afunc(\rflow)},
$$
whence \ $\afunc(\vars) \,\cdot \,\rflow'_{\vars}(\vars,t) \, = \, \afunc\circ\rflow(\vars,t)$.
\endproof

On the other hand, 
$$\afunc\circ\rflow(\vars,t)\, = \, \afunc \, + \, (\rflow-\vars)\, \afunc'_{\vars} \, + \, (\rflow-\vars)^2 \, \mfunc(\vars,t),$$
where $\mfunc(\vars,t)$ is a certain smooth function.

Hence~\eqref{equ:cfunc_main_rel} can be rewritten in the following form:
$$
\cfunc\,'_{\vars} = \frac{\afunc\circ\rflow - \afunc - (\rflow-\vars) \afunc\,'_{\vars}}{t\,\afunc^2}=
\frac{ (\rflow-\vars)^2 }{t^2\,\afunc^2} \, t \mfunc(\vars,t)=
\cfunc^2\, t \mfunc(\vars,t),
$$
which gives us~\eqref{equ:dif_eq_cfunc}.
Notice that $\cfunc(\vars,0)=\frac{1}{c(0)}$, whence 
$\cfunc(\vars,t) = \frac{1}{c(0)} + t \bfunc(\vars,t)$ for some smooth function $\bfunc(\vars,t)$.
Therefore
\begin{equation}\label{lm:char_flowR2_details}
\rflow(\vars,t) \, = \, \vars \, + \, \frac{t \, \afunc(\vars)}{c(0)} \, + \, t^2 \, \afunc(\vars) \, \bfunc(\vars,t).
\end{equation}

Comparing this with~\eqref{equ:flow_decomp} we obtain that $c(0)=1$.

Lemma~\ref{lm:cfunc_rel} and Proposition~\ref{pr:char_flowR1} are proved.
\endproof

\section{Stabilizers of local left-right action}\label{sect:local_act}
Let $\smzrm$ be the algebra of germs of smooth functions $\Rm\to\RRR$ at $0\in\RRR^{m}$ and $\mmmRm$ a unique maximal ideal of $\smzrm$ consisting of functions $\mrsfunc$ such that $\mrsfunc(0)=0$.
Let also $\DiffRmz$ be the groups of germs of preserving orientation diffeomorphisms $\difM:\Rm\to\Rm$ at $0\in\Rm$ such that $\difM(0)=0$.

Designate respectively by $\smzr$, $\mmmR$, and $\DiffRz$ the analogous objects for $\RRR$.
Then the groups $\DiffRmz$ and $\DiffRmz\times\DiffRz$ acts on $\mmm$ by formulas~\eqref{equ:act_M} and ~\eqref{equ:act_MP} respectively.  
We will call these actions {\em local right\/} and {\em local left-right\/} respectively.

We also say that $\mrsfunc,\gfunc\in\smzrm$ are {\em left-right\/} ({\em right}) equivalent iff they belong to the same orbit with respect to a local left-right (right) action.
 
For $\mrsfunc\in\mmmRm$ let 
$$
\StabfRmR=\{(\difM,\difR)\in\DiffRmz\times\DiffRz \ | \ \difR\circ\mrsfunc=\mrsfunc\circ\difM \}
$$ be the stabilizer of $\mrsfunc$ with respect to this action.

\subsection{Property \condJacIda{\afunc}.} 
Notice that each $\mrsfunc\in\mmmRm$ yields a homomorphism of algebras:
$$\mrsfunc^{*}:\smzr\to\smzrm, \qquad \mrsfunc^{*}(\afunc)=\afunc\circ\mrsfunc.$$

Denote by $\Jacidfzero\subset\smzrm$ the Jacobi ideal of $\mrsfunc$ at $0\in\Rm$, i.e.\!\! the ideal generated by partial derivatives $\mrsfunc^{\prime}_{x_1}, \ldots, \mrsfunc^{\prime}_{x_{\dimM}}$  of $\mrsfunc$.

\begin{defn}\label{defn:cond_Ja}
Let $\afunc\in\mmmR$. 
We say that $\mrsfunc\in\mmmRm$ has property \condJacIda{\afunc}\ at $0\in\Rm$ if 
\begin{equation}\label{equ:cond_Ja} 
\mrsfunc^{*}(\afunc) \ = \ \afunc\circ\mrsfunc \ \in \ \Jacidfzero.
\end{equation}
Equivalently, there exists a vector field $\Fld$ at $0\in\Rm$ such that 
$$ \dfval{\mrsfunc}{\Fld} = \afunc\circ\mrsfunc: \Rm \stackrel{\mrsfunc}{\to} \RRR \stackrel{\afunc}{\to} \RRR.$$

For instance, if $\afunc=\id_{\RRR}$, then \condJacIda{\id}\ means that $\mrsfunc\in\Jacidfzero$, which is precisely the condition \condJacId. 
More generally, if $\afunc(\vars)=\vars^{k}$, then \condJacIda{\vars^{k}} means that $\mrsfunc^{k}\in\Jacidfzero$.

Notice also that for every $\bfunc\in\smzr$ the condition \condJacIda{\afunc}\ implies \condJacIda{\afunc\bfunc}, i.e.\!\! if $\afunc\circ\mrsfunc=\dfval{\mrsfunc}{\Fld}$, then 
$(\afunc\bfunc)\circ\mrsfunc=(\afunc\circ\mrsfunc)(\bfunc\circ\mrsfunc) = \dfval{\mrsfunc}{(\bfunc\circ\mrsfunc\cdot \Fld)}\in\Jacidfzero$. 
Hence, if $\bfunc(0)\not=0$, then \condJacIda{\afunc}\ and \condJacIda{\afunc\bfunc} are equivalent.
\end{defn}

\begin{lem}\label{lm:condJac_inv_orb}
{\rm(1)}~Property \condJacIda{\afunc}\ is invariant with respect to local {\em right\/} equivalence.

{\rm(2)}~Property \condJacIda{\vars^{k}}\ is invariant with respect to local {\em left-right\/} equivalence.
\end{lem}
\proof
(1) Suppose that $\mrsfunc$ has \condJacIda{\afunc}, i.e. $\afunc\circ\mrsfunc = \dfval{\mrsfunc}{\Fld}$ for some vector field $\Fld=(\Fld_1,\ldots,\Fld_{\dimM})$ at $0\in\Rm$.
We have to prove that for each $\difM\in\DiffRmz$ the function $\gfunc=\mrsfunc\circ\difM$ has \condJacIda{\afunc} as well.
Indeed, since $\nabla\gfunc=\nabla(\mrsfunc\circ\difM)=Th\cdot(\nabla f)\circ h$ we get
$$
\afunc\circ\gfunc = \afunc\circ\mrsfunc\circ\difM=\dfvalx{\mrsfunc}{\Fld}{\difM} =
\sum_{i=1}^{\dimM}  (\Fld_i \circ\difM) \cdot (\mrsfunc'_{x_i}\circ\difM) =
\dfval{\gfunc}{\left[ (\Fld \circ\difM) \cdot (T\difM)^{-1} \right]} \ \in \ \JFzero{\gfunc}.
$$

(2) Suppose that $\mrsfunc$ has \condJacIda{\vars^{k}}, i.e. $\mrsfunc^{k}=\dfval{\mrsfunc}{\Fld}$ for some vector field $\Fld$ at $0\in\Rm$.
We have to prove that for each $(\difR,\difM)\in\DiffRmz\times\DiffRz$ the function $\gfunc=\difR\circ\mrsfunc\circ\difM^{-1}$ has \condJacIda{\vars^{k}}.
Due to (1) we can assume that $\difM=\id_{\Rm}$, whence $\gfunc=\difR\circ\mrsfunc$.
Then
$$
\dfval{\gfunc}{\Fld} = \sum_{i=1}^{\dimM} \Fld_i\cdot(\difR\circ\mrsfunc)'_{x_i}  =
\difR'(\mrsfunc) \ \sum_{i=1}^{\dimM} \Fld_i\,\mrsfunc'_{x_i}  =
\difR'(\mrsfunc) \cdot \dfval{\mrsfunc}{\Fld} = 
\difR'(\mrsfunc) \, \mrsfunc^{k}.
$$
Notice that $\difR(\vars)=\bdifR(\vars)\vars$, where $\bdifR\in\smzr$ and $\bdifR(0)=\difR'(0)>0$. Therefore
$$ 
\gfunc^{k} = (\difR\circ\mrsfunc)^{k}=
\bdifR(\mrsfunc)^{k} \, \mrsfunc^{k} =
 \frac{\bdifR(\mrsfunc)^{k}}{\difR'(\mrsfunc)} \ \difR'(\mrsfunc) \mrsfunc^{k} =
 \dfval{(\difR\circ\mrsfunc)}{ \left( \frac{\bdifR(\mrsfunc)^{k}}{\difR'(\mrsfunc)} \ \Fld \right)} \ \in \ \JFzero{\gfunc}. \qed
$$

\begin{rem}
(c.f.\!\! Corollary~\ref{cor:Vdf_ska_is_normal_in_Vdf_a}).
It seems that for arbitrary $\afunc$ the property \condJacIda{\afunc}\ is not a local left-right invariant.
The crucial property of the function $\afunc(\vars)=\vars^{k}$ that we have used for the proof of (2)  is that 
$\afunc(\difR\circ\mrsfunc) = \omega \afunc(\mrsfunc)$, i.e.\!\! $(\difR\circ\mrsfunc)^{k}=\omega \mrsfunc^{k}$ for some $\omega\in\smzr$.
This property may not hold in general.
For example, let 
$$
\afunc(\vars)=\left\{
\begin{array}{ll}
e^{-1/\vars}, & \vars>0, \\
0, & \vars\leq 0,
\end{array}
\right.
\qquad
\difR(\vars)=2\vars,
\qquad
\mrsfunc(\vars)=\vars.
$$
Then 
$\afunc\circ\difR\circ\mrsfunc(\vars)=e^{-1/2\vars}=e^{1/2\vars} \cdot \afunc\circ\mrsfunc(\vars)$ for $\vars>0$ and the function $\omega(\vars)=e^{1/2\vars}$ does not extend to a smooth function near $0$.
\end{rem}

\subsubsection{}
The following theorem establishes the relationships between the group $\Vdf{\afunc}$, the condition \condJacIda{\afunc} for $\mrsfunc$, and the stabilizer $\StabfRmR$. 
 
Let $\prMPP:\DiffRmz\times\DiffRz\to\DiffRz$ be natural projection.

\begin{thm}\label{th:ex_loc_sect}
Suppose that $\mrsfunc\in\mmmRm$ has property \condJacIda{\afunc}\ for some admissible $\afunc\in\mmmR$. 
Then $\Vdf{\afunc} \subset \prMPP(\StabfRmR)$ and $\prMPP$ has a section $\sectMPP:\Vdf{\afunc} \to \StabfRmR$ being a homomorphism.

Equivalently, if \ $\dfval{\mrsfunc}{\mFld} = \afunc\circ\mrsfunc$ \ for some vector field $\mFld$ at $0\in\Rm$, then there exists a homomorphism $\difShift:\Vdf{\afunc}\to\DiffRmz$ such that $$\difR\circ\mrsfunc=\mrsfunc\circ\difShift(\difR),$$
i.e. $(\difShift(\difR),\difR)\in\StabfRmR$.
Thus each $\difM\in\Vdf{\afunc}$ is \LT-trivial for $\mrsfunc$.

Moreover, let $\mflow$ be a flow generated by $\mFld$. 
Then $\difShift$ is defined by the following formula:
$\difShift(\difR)(x) = \mflow(x,\shfunc{\difR}\circ\mrsfunc(x))$, where $\difR\in\Vdf{\afunc}$ and $\shfunc{\difR}\in\smzr$ is a shift-function of $\difR$ with respect to the flow generated by the vector field $\rFld(\vars)=\afunc(\vars)\frac{d}{ds}$.
\end{thm}

As a particular case of this theorem we obtain the following local variant of Theorem~\ref{th:Stab}.
\begin{thm}\label{th:loc_Stab}
Suppose that $\mrsfunc\in\mmmRm$ satisfies \condJacIda{\id}, i.e.\!\! $\mrsfunc\in\Jacidfzero$.
Then $\prMPP(\StabfRmR)=\DiffRz=\Vdf{\id}$ and $\prMPP$ admits a section 
$\sectMPP:\DiffRz \to \StabfRmR$ being a homomorphism.
Thus each $\adif\in\DiffRz$ is \LT-trivial for $\mrsfunc$.
\end{thm}

\begin{cor}
Suppose that $\mrsfunc,\gfunc\in\minid{\Rm}$ satisfy \condJacIda{\id}.
Then $\mrsfunc$ and $\gfunc$ are left-right equivalent iff they are right equivalent.
\end{cor}
\proof
If $\mrsfunc$ and $\gfunc$ are left-right equivalent, i.e. 
$\gfunc = \difR\circ\mrsfunc\circ\difM^{-1}$ for a certain $(\difM,\difR)\in\DiffRmz\times\DiffRz$, 
then $\gfunc = \mrsfunc\circ\difShift(\difR)\circ\difM^{-1}$.
The converse statement is evident.
\endproof

We shall prove Theorem~\ref{th:ex_loc_sect} in Section~\ref{sect:prf:th:ex_loc_sect}.

\subsection{Characterization of property \condJacIda{\afunc}.}
Let $\mFld$ be a vector field on $\Rm$ and $\mflow:\mnbh\times\Interv\to\Rm$ a local flow generated by $\mFld$, where $\mnbh$ is a neighborhood of $0\in\Rm$ and $\Interv$ is an open interval containing $0\in\RRR$.

Let $\afunc:\RRR\to\RRR$ be a smooth function, 
$\rFld(\vars)=\afunc(\vars)\frac{d}{d\vars}$ a vector field on $\RRR$, and 
$\rflow:\rnbh\times\Interv\to\RRR$ a local flow generated by $\rFld$, where $\rnbh$ is a neighborhood of $0\in\RRR$.
Notice that $\rFld$ is a section of a (trivial) tangent bundle $T\RRR$ defined by 
$$ \rFld:\RRR \to T\RRR \equiv \RRR\times\RRR, \qquad \rFld(\vars) = (\vars,\afunc(\vars)).$$

\begin{lem}\label{lm:dfG_af}
Suppose that $\mrsfunc:\Rm\to\RRR$ is a smooth function such that $\mrsfunc(0)=0$.
Then the following conditions on $\mFld$, $\rFld$, $\mflow$, $\rflow$, and $\afunc$ are equivalent:

\begin{enumerate}
\item
The following diagram is commutative:
$$
\begin{array}{ccccc}
\begin{CD}
T\mnbh @>{T\mrsfunc}>> T\rnbh \\ 
@A{\mFld}AA @AA{\rFld}A \\
\mnbh @>{\mrsfunc}>> \rnbh
\end{CD}
& \qquad &
\text{i.e.}
& \qquad &
 \dfval{\mrsfunc}{\mFld}=\afunc\circ\mrsfunc,
\end{array}
$$
which is precisely condition \condJacIda{\afunc}.
\item
For each $t\in\Interv$ the following diagram is commutative: 
$$
\begin{array}{ccccc}
\begin{CD}
\mnbh @>{\mrsfunc}>> \rnbh \\
@V{\mflow_t}VV @VV{\rflow_t}V \\
\Rm @>{\mrsfunc}>> \RRR
\end{CD}
& \qquad &
\text{i.e.}
& \qquad &
\rflow(\mrsfunc(x),t)=\mrsfunc\circ\mflow(x,t);
\end{array}
$$

\item
For each smooth function $\sfunc:\rnbh\to\Interv$ define two mappings
$$
\begin{array}{lcl}
\difM_{\sfunc}:\mnbh\to\Rm & & \difM(x)=\mflow(x,\sfunc\circ\mrsfunc(x)),  \\ [1.3mm]
\difR_{\sfunc}:\rnbh\to\RRR & & \difR(\vars)=\rflow(\vars,\sfunc(\vars)).
\end{array}
$$
Then the following diagram is commutative:
$$
\begin{array}{ccccc}
\begin{CD}
\mnbh @>{\mrsfunc}>> \rnbh \\
@V{\difM_{\sfunc}}VV @VV{\difR_{\sfunc}}V \\
\Rm @>{\mrsfunc}>> \RRR
\end{CD}
& \qquad &
\text{i.e.}
& \qquad &
\difR_{\sfunc}\circ\mrsfunc=\mrsfunc\circ\difM_{\sfunc}.
\end{array}
$$
In this case $\difM_{\sfunc}$ is an embedding iff $\difR_{\sfunc}$ is.

\item
There exist smooth isotopies 
$$
\bar\mflow:\mnbh\times[0,1]\to\Rm \qquad \text{and} \qquad \bar\rflow:\rnbh\times [0,1]\to\RRR
$$
such that 
$\bar\rflow_0=\id_{\rnbh}$, $\bar\mflow_0=\id_{\mnbh}$, 
$\bar\rflow_t\circ\mrsfunc=\mrsfunc\circ\bar\mflow_t$,
\begin{equation}\label{equ:cond4}
\frac{\partial \bar\rflow}{\partial t}(\vars,0)=\afunc(\vars), \qquad \text{and} \qquad 
\frac{\partial \bar\mflow}{\partial t}(\vars,0)=\mFld(\vars).
\end{equation}
\end{enumerate}
\end{lem}
\proof
(1)$\Rightarrow$(2)
Let $x\in\mnbh$.
It suffices to prove that if $\omega(t)=\rflow(x,t)$ is a trajectory of $\mFld$, i.e.
$\omega'_{t}(t)=\mFld(\omega(t))$, then $\mrsfunc\circ\omega(t)$ is a trajectory of $\rFld$, i.e. $(\mrsfunc\circ\omega(t))'_{t}=\rFld(\mrsfunc\circ\omega(t))$.
Indeed, 
$$(\mrsfunc\circ\omega(t))'_{t} \,=\,
\dfval{\mrsfunc}{\omega'_{t}(t)} \,=\, 
\dfvalx{\mrsfunc}{\mFld}{\omega(t)} \,\stackrel{(1)}{=}\, 
\rFld(\mrsfunc\circ\omega(t)).$$

(2)$\Leftrightarrow$(3)
Statement (3) can be obtained by substituting into (2) the function $\sfunc\circ\mrsfunc(x)$ instead of $t$.
Conversely, (2) is a particular case of (3) for the constant function $\sfunc(\vars)=t$.

It remains to prove that $\difM_{\sfunc}$ is an embedding iff $\difR_{\sfunc}$ is. By~\cite[Theorem~19]{Maks:Shifts} $\difM_{\sfunc}$ (resp. $\difR_{\sfunc}$) is an embedding that preserves orientations of trajectories of $\rflow$ iff $\dfval{(\sfunc\circ\mrsfunc)}{\mFld}>-1$ (resp. $\dfval{\sfunc}{\rFld}>-1$).
But under assumption (1) these expressions coincide:
$$
\dfval{(\sfunc\circ\mrsfunc)}{\mFld} \,=\,
\sfunc'(\mrsfunc)\,\cdot\,\dfval{\mrsfunc}{\mFld} \,=\,
\sfunc'(\mrsfunc)\,\cdot\,\rFld\circ\mrsfunc \,=\,
\dfvalx{\sfunc}{\rFld}{\mrsfunc}. 
$$

(2)$\Rightarrow$(4)
It suffices to set $\bar\rflow=\rflow$ and $\bar\mflow=\mflow$, where $\rflow$ and $\mflow$ are the corresponding flows. 
Then condition~\eqref{equ:cond4} simply means that $\rflow$ and $\mflow$ are generated by $\rFld$ and $\mFld$ respectively.

(4)$\Rightarrow$(1)
Let $\mFld_{t}(x) = \bar\mflow'_{t}(x,t)$ 
be the one parametric family of vector fields on $\mnbh$ corresponding to the isotopy $\bar\mflow_t$. 
Differentiating \ $\mrsfunc\circ\bar\mflow(x,t) = \bar\rflow(\mrsfunc(x),t)$ \ in $t$ we obtain
$$
\dfvalx{\mrsfunc}{\mFld_{t}}{\bar\mflow(x,t)} = \bar\rflow'_{t}(\mrsfunc(x),t).
$$
For $t=0$ we have that $\bar\mflow(x,0)=x$ and $\bar\rflow'_{t}(\vars,t)=\afunc(\vars)$, therefore 
$\dfval{\mrsfunc}{\mFld_0}=\afunc\circ\mrsfunc$.
\endproof

Let us mention one particular case of Lemma~\ref{lm:dfG_af} which will play a crucial role in the proof of Theorem~\ref{th:Stab}.
For the convenience of the reader we formulate this statement explicitely.
\begin{lem}\label{lm:f_flow__f_exp}
Let $\mrsfunc\in\mmmRm$, $\eps\in\RRR$, $\afunc(\vars)=\eps\vars$, 
$\rFld(\vars)=\eps\vars\frac{d}{ds}$ be a vector field on $\RRR$,
and $\rflow(\vars,t)=\vars e^{\eps t}$ the flow generated by $\rFld$.
Then the following conditions are equivalent, as they are the conditions (1) and (2) of Lemma~\ref{lm:dfG_af} for this case:
\begin{gather}
  \dfval{\mrsfunc}{\mFld} = \eps \mrsfunc, \label{equ:j_main_relation} \\
  \mrsfunc \circ \mflow(x,t) = \mrsfunc(x) \cdot e^{\eps t}\label{equ:main_relation}. \end{gather}
\end{lem}
 
\subsection{Proof of Theorem~\ref{th:ex_loc_sect}}
\label{sect:prf:th:ex_loc_sect}
Suppose that $\dfval{\mrsfunc}{\Fld}=\afunc\circ\mrsfunc$, where $\mFld$ is a vector field on $\Rm$.
Let $\mflow$ be a local flow on $\Rm$ at $0\in\Rm$ generated by $\mFld$, and 
$\rflow$ be a local flow on $\RRR$ at $0\in\RRR$ generated by the vector field $\rFld(\vars)=\afunc(\vars)\frac{d}{ds}$.

If $\difR\in\Vdf{\afunc}$, then by Theorem~\ref{th:Vdf_shifts} $\difR(\vars)=\rflow(\vars,\sfunc(\vars))$ for some smooth function $\sfunc\in\smzr$.
Define $\difShift(\difR)$ by $\difShift(\difR)(x)=\mflow(x,\sfunc\circ\mrsfunc(x))$.
Then by (3) of Lemma~\ref{lm:dfG_af} we have $\difR\circ\mrsfunc=\mrsfunc\circ\difShift(\difR)$.

It remains to show that $\difShift$ is a homomorphism.
Let $\difR_i(\vars)=\rflow(\vars,\sfunc_i(\vars))$ and $\difM_i(x)=\mflow(x,\sfunc_i\circ\mrsfunc(x))$.\
Then 
$$
\begin{array}{lcl}
\difR_2\circ\difR_1(\vars) & = & \rflow( \difR_1(\vars), \sfunc_2\circ\difR_1(\vars)) =
\rflow( \vars, \underbrace{\sfunc_1(\vars) + \sfunc_2\circ\difR_1(\vars)}_{\sfunc(\vars)}), \\
\difM_2\circ\difM_1(x) & = & 
\mflow( \difM_1(x), \sfunc_2\circ\mrsfunc\circ\difM_1(x)) = \\
& & 
\mflow( x, \sfunc_1\circ\mrsfunc(x) + \sfunc_2\circ\difR_1\circ\mrsfunc(x)) = 
\mflow(x,\sfunc\circ\mrsfunc(x)). \qed
\end{array}
$$

\subsection{Problems.}
By Theorem~\ref{th:ex_loc_sect} we know that if $\mrsfunc$ satisfies \condJacIda{\afunc}, then $\Vdf{\afunc}$ is included in the group $\prMPP(\StabfRmR)$ of all \LT-trivial diffeomorphisms for $\mrsfunc$.
In the next section we will prove that a very large class of functions $\mrsfunc$ satisfies the ``maximal'' condition \condJacIda{\id} which is the same as \condJacId.
The property \condJacId\ appears typical since it holds for non-degenerate critical points, all simple singularities $A_k$, $D_k$, $E_6$, $E_7$, $E_8$, and even for formal series.
Thus for these functions by  Theorem~\ref{th:ex_loc_sect} we have that $\prMPP(\StabfRmR)=\DiffRz=\Vdf{\id}$.
Moreover, we also show that \condJacId\ is invariant with respect to the stable equivalence of singularities (Definition~\ref{defn:stable_equiv}).

On the other hand, there are singularities that do not satisfy \condJacId, see Claim~\ref{clm:ex_not_in_J}.
This yields the following problems whose solutions would probably give new invariants of pathological singularities and in particular of flat functions.

1) {\em
If $\Vdf{\afunc}\subset\prMPP(\StabfRmR)$, does the condition \condJacIda{\afunc}\ holds for $\mrsfunc$?
}

2) Suppose that $\Vdf{\vars\afunc} \subseteq \prMPP(\StabfRmR) \subseteq \Vdf{\afunc}$.
Recall that $\Vdf{\afunc}/\Vdf{\vars\afunc}\approx\RRR$, thus the group $G=\prMPP(\StabfRmR)/\Vdf{\vars\afunc}$ is a subgroup in $\RRR$.
{\em Is $G$ closed or everywhere dense in $\RRR$? Does it coincide with either $0$ or $\RRR$?}
In the last two cases $\prMPP(\StabfRmR)$ coincides with either $\Vdf{\vars\afunc}$ or $\Vdf{\afunc}$.

3) {\em 
Is it true in general, that for every $\mrsfunc\in\mmmRm$ we have that $\prMPP(\StabfRmR)$ coincides with some group $\Vdf{\afunc}$?
}

\section{Condition \condJacIda{\id}}\label{sect:cond_JacId}
In this section we give examples of singularities having property \condJacIda{\id}.
This property is invariant with respect to a stable equivalence (Corollary~\ref{cor:Jid_inv_stable_eq}) and holds for non-degenerate, simple singularities (Corollary~\ref{cor:Jid_simple_sing}) and also for formal series.
We also show that there are singularities for which \condJacIda{\id} fails.

Let $\smzrm$ be the algebra of germs of smooths functions $\Rm\to\RRR$ at $0\in\Rm$, $\mmmRm$ a unique maximal ideal of $\smzrm$ consisting of functions $\mrsfunc$ such that $\mrsfunc(0)=0$, and $\Jacidfzero$ the Jacobi ideal of $\mrsfunc$ generated by partial derivatives of $\mrsfunc$.

Recall that $\mrsfunc$ has property \condJacIda{\id}\ at $0\in\Rm$ provided $\mrsfunc\in\Jacidfzero$.
Moreover, by (2) of Lemma~\ref{lm:condJac_inv_orb} \  \condJacIda{\id}\ is the property of  the {\em local left-right orbit\/} of $\mrsfunc$.

\begin{lem}\label{lm:class_JM_dim1}
Let $\dimM=1$. Then $\mrsfunc$ has \condJacIda{\id}\ if and only if
$\mrsfunc(x) = \afunc(x)\mrsfunc'(x)$ for some $\afunc\in\smzr$. \qed
\end{lem}

\begin{lem}\label{lm:examples_func_in_jac}
Suppose that $\gfunc\in\mmmRm$ has \condJacIda{\id}, so $\dfval{\gfunc}{\Fld}=\gfunc$ for some vector field $\Fld$ at $0$.
Let also $\mrsfunc\in\mmmRm$ satisfies one of the following conditions:
\begin{enumerate}
\item
$\mrsfunc = \gfunc^{a}$ \ for some \ $a\in\RRR\setminus0$ (emphasize that $\mrsfunc$ is assumed $C^{\infty}$);
\item
$\mrsfunc(x)=\left\{
\begin{array}{ll}
e^{-\frac{1}{|\gfunc(x)|}}, & \gfunc(x)\not=0 \\
0, & \gfunc(x)=0;
\end{array}
\right.
$

\end{enumerate}
Then $\mrsfunc$ also has property \condJacIda{\id}, i.e. $\mrsfunc\in\Jacidgzero$.
\end{lem}
\proof
We will show that in all the cases $\dfval{\mrsfunc}{\bfunc\Fld}=\mrsfunc$ for some $\bfunc\in\smzr$.

\smallskip

(1)
$\dfval{\mrsfunc}{\frac{1}{a}\Fld} =
\dfval{(\gfunc^{a})}{\frac{1}{a}\Fld} = \gfunc^{a-1} \dfval{\gfunc}{\Fld} =
\gfunc^{a-1}\gfunc = \mrsfunc$.

(2)
It is well known that the following function
$\psi(x)=\left\{
\begin{array}{ll}
e^{-\frac{1}{|x|}}, & x\not=0 \\
0, & x=0
\end{array}
\right.$
is smooth, whence so is $\mrsfunc = \psi\circ\gfunc$.
Then
$$
\dfval{\mrsfunc}{\gfunc\Fld} =
\frac{\dfval{\gfunc}{\gfunc\Fld}}{\gfunc^{2}} e^{-\frac{1}{|\gfunc|}} =
e^{-\frac{1}{|\gfunc|}} = \mrsfunc. \ \qed
$$

\begin{lem}\label{lm:class_JM}
Let $\mrsfunc\in\mmmRm$.
Suppose that in some local coordinates at $0$ the function $\mrsfunc$ satisfies one of the following conditions:
\begin{enumerate}
\item[{\rm(i)}]
$\mrsfunc(x_1,\ldots,x_{\dimM})=x_1$, which is equivalent to the assumption that $0$ is a regular point of $\mrsfunc$;

\item[{\rm(ii)}]
$\mrsfunc(t x_1,\ldots,t x_{\dimM})=t^{n}\mrsfunc(x_1,\ldots,x_{\dimM})$
for all $t>0$, i.e.\!\! $\mrsfunc$ is homogeneous;

\item[{\rm(iii)}]
$\mrsfunc(x_1,\ldots,x_{\dimM})=
x_1^{a_1} \pm x_{2}^{a_2} \pm \cdots \pm x_{k}^{a_{k}}$; 
\medskip

\item[{\rm(iv)}]
$\mrsfunc(x_1,\ldots,x_{\dimM}) =
x_1^{a_1} + x_{1}^{b_1} x_{2}^{a_2} + x_{2}^{b_2} x_{3}^{a_3} + \cdots +
 x_{k-1}^{b_{k-1}} x_{k}^{a_{k}}$;

\medskip
\item[{\rm(v)}]
$m=1$ and $\mrsfunc^{(s)}(0)\not=0$ for some $s\geq 1$, i.e. $\mrsfunc$ is not flat at $0$;
\end{enumerate}
where $a_i\geq1$, $b_j\geq0$, and $k\leq \dimM$.
Then $\mrsfunc$ has property \condJacIda{\id}.
\end{lem}
\proof
{\rm(i)} $\mrsfunc= x_1 \mrsfunc'_{x_1}$.

\smallskip

{\rm(ii)}
This case follows from the well-known Euler identity for homogeneous functions:
$ \mrsfunc(x) = \frac{1}{n} \mathop\sum\limits_{i=1}^{\dimM} x_i \mrsfunc'_{x_{i}}.$
Let us recall its proof.
Define a vector field by the formula
$\Fld(x_1,\ldots,x_{\dimM}) = \frac{1}{n}(x_1,\ldots,x_{\dimM}).$
Then
\begin{multline}\label{equ:Euler_homog_id}
\dfvalx{\mrsfunc}{\Fld}{x} =
\lim\limits_{s\to 0} \frac{\mrsfunc(x+s\Fld(x)) - \mrsfunc(x)}{s} = \\ =
\lim\limits_{s\to 0} \frac{\mrsfunc((1+s/n)x) - \mrsfunc(x)}{s} =
\lim\limits_{s\to 0} \frac{(1+s/n)^n - 1}{s} \mrsfunc(x)  = \mrsfunc(x).
\end{multline}

\medskip

{\rm(iii)}
$\mrsfunc = \frac{x_1}{a_1} \mrsfunc'_{x_1} \pm
\frac{x_2}{a_2} \mrsfunc'_{x_2} \pm \cdots \pm
\frac{x_{k}}{a_{k}} \mrsfunc'_{x_{k}}$.

\medskip

{\rm(iv)}
$\mrsfunc = \frac{x_1}{a_1} \mrsfunc'_{x_1} + \left( 1-\frac{b_1}{a_1} \right) \frac{x_2}{a_2} \mrsfunc'_{x_2} + \left[ 1- \frac{b_2}{a_2}\left( 1-\frac{b_1}{a_1} \right) \right] \frac{x_{3}}{a_{3}}  \mrsfunc'_{x_{3}} + \ldots$

\medskip

{\rm(v)}
Let $s\geq 1$ be the first number for which $\mrsfunc^{(s)}(0)\not=0$.
Then $\mrsfunc$ is right equivalent to $x^n$, and our statement follows from (ii).
\endproof
\begin{cor}\label{cor:Jid_simple_sing}
Non-degenerate singularities $\sum \pm x_i^2$ and simple singularities $A_k(x)=x^k$ $(k\geq1)$, $D_k(x,y)=x^2y+y^{k-1}$ $(k\geq4)$, $E_6(x,y)=x^3+y^4$, $E_7(x,y)=x^3+xy^3$, $E_8(x,y)=x^3+y^5$, see e.g.~\cite{AVG},
have \condJacIda{\id}.\qed
\end{cor}

\begin{lem}\label{lm:stable_equiv}
Let $\mrsfunc\in\minid{\Rm}$, $\gfunc\in\minid{\Rn}$. 
Define $\hfunc\in\minid{\RRR^{m+n}}$ by the formula 
$\hfunc(x,y)=\mrsfunc(x)+\gfunc(y)$,\ where $x=(x_1,\ldots,x_{m})\in\Rm$ and $y=(y_1,\ldots,y_{n})\in\Rn$.
Then $\mrsfunc$ and $\gfunc$ have property \condJacIda{\id}\ iff $\hfunc$ has.
\end{lem}
\proof
Suppose that $\mrsfunc(x)=\dfvalx{\mrsfunc}{\fFld}{x}$ and $\gfunc(y)=\dfvalx{\gfunc}{\gFld}{y}$ for certain vector fields $\fFld$ and $\gFld$ defined on $\Rm$ and $\Rn$ respectively.
We can regard these vector fields as components of the following vector field 
$\hFld(x,y)=(\fFld(x),\gFld(y))$ on $\RRR^{m+n}$.
Then $\dfval{\mrsfunc}{\gFld} = \dfval{\gfunc}{\fFld}=0$, whence
$$ \dfvalx{\hfunc}{\hFld}{x,y} =
\dfvalx{\mrsfunc}{\fFld}{x} + \dfvalx{\gfunc}{\gFld}{y} =
\mrsfunc(x) + \gfunc(y) = \hfunc(x,y).
$$

Conversely, suppose that $\hfunc(x,y)=\dfval{\hfunc}{\hFld}(x,y)$, where
$$\hFld(x,y)=(\fFld(x,y),\gFld(x,y))$$ is a certain vector field on $\RRR^{m+n}$.
Notice that 
$$
\dfval{\hfunc}{\hFld}(x,y) = 
\dfvalx{(\mrsfunc+\gfunc)}{(\fFld,\gFld)}{x,y}=
\sum_{i=1}^{m} \mrsfunc'_{x_i}(x) \fFld_i(x,y) +
\sum_{j=1}^{n} \gfunc'_{y_j}(y) \gFld_j(x,y),
$$
where $\fFld=(\fFld_1,\ldots,\fFld_{m})$ and $\gFld=(\gFld_1,\ldots,\gFld_{n})$.

Let $\bar\fFld(x)=\fFld(x,0)$ and $\bar\gFld(y)=\gFld(0,y)$ be vector fields on $\Rm$ and $\Rn$ respectively.
By (i) of Lemma~\ref{lm:class_JM} we can assume that $0\in\Rm$ and $0\in\Rn$ are critical points for $\mrsfunc$ and $\gfunc$ respectively.
Then $\mrsfunc'_{x_i}(0)=0$ and $\gfunc'_{y_j}(0)=0$, whence 
$$
\mrsfunc(x) = \hfunc(x,0) = \dfval{\hfunc}{\hFld}(x,0) =
\sum_{i=1}^{m} \mrsfunc'_{x_i}(x) \fFld_i(x,0) = \dfvalx{\mrsfunc}{\bar\fFld}{x}
$$
and similarly $\gfunc(y) = \dfvalx{\gfunc}{\bar\gFld}{y}$.
\endproof

\begin{defn}[see e.g.~\cite{AVG}]\label{defn:stable_equiv}
Two functions $\mrsfunc_{k}\in\minid{\RRR^{m_k}}$ $k=1,2$, are {\em stably equivalent\/} if there are non-degenerated quadratic forms 
$$\gfunc_{k}(y_1,\ldots,y_{n_{k}})=\mathop\sum\limits_{i=1}^{n_k}\pm y_i^{2}, \qquad k=1,2$$ such that 
$m_1+n_1=m_2+n_2$ and functions 
$\mrsfunc_1+\gfunc_1$ and $\mrsfunc_2+\gfunc_2$ defined on $\RRR^{m_1+n_1}$ are {\em right\/} equivalent (belong to the same orbit with respect to local {\em right\/} action).
\end{defn}
By either (ii) or (iii) of Lemma~\ref{lm:class_JM} each quadratic form has property \condJacIda{\id}, whence from Lemma~\ref{lm:stable_equiv} we get the following statement:
\begin{cor}\label{cor:Jid_inv_stable_eq}
Suppose that functions $\mrsfunc_{k}\in\minid{\RRR^{m_k}}$ $(k=1,2)$ are stably equivalent.
Then $\mrsfunc_1$ has property \condJacIda{\id}\ iff $\mrsfunc_2$ has.
\qed
\end{cor}

\begin{prop}
Let 
$$
\mrsfunc = 
\sum\limits_{i_1,\ldots, i_{\dimM} \geq 0} 
a_{i_1 \ldots i_{\dimM}} \, x_1^{i_1} \ldots x_{\dimM}^{i_{\dimM}}
\ \in \ \RRR[[x_1,\ldots,x_{\dimM}]]
$$
 be a formal series without first term, i.e. $\mrsfunc(0)=0$.
Then $\mrsfunc$ has \condJacIda{\id}\ in $\RRR[[x_1,\ldots,x_{\dimM}]]$.
\end{prop}
\proof
We should find formal series $\Fld_1,\ldots,\Fld_{\dimM}\in\RRR[[x_1,\ldots,x_{\dimM}]]$ such that $$\mrsfunc=\sum_{i=1}^{\dimM}\mrsfunc'_{x_i} \Fld_i.$$
This relation gives a system of linear equations on the coefficients of $\Fld_i$ such that the solution can be found recurrently.
The details are left to the reader.
\endproof

This statement gives a hope, that \condJacIda{\id}\ holds for analytical functions, and in particular for polynomials (as usual, we have to seek for converging series $\Fld_i$), but I do not know is it true.
On the other hand, the following statement shows that \condJacIda{\id}\ can fail in the non-analytical case.

\begin{claim}\label{clm:ex_not_in_J}
Let $\mrsfunc\in\mmmRm$.
Suppose that there is a sequence $\{\pnt_{i}\}$ of critical points of $\mrsfunc$ converging to $0$ and such that $\mrsfunc(\pnt_{i})\not=0$, in particular, the critical {\em value\/} $0$ of $\mrsfunc$ is not isolated.
Then $\mrsfunc$ does not satisfy \condJacIda{\id}, i.e. \!\! $\mrsfunc\not\in\Jacidfzero$.
\end{claim}
\begin{rem}
Such a function exists. For example let
$$
\gfunc(x)=\left\{
\begin{array}{ll}
e^{-\frac{1}{\sin^2 {1/x}}}, & x \not=\frac{1}{\pi n}, \\ [2mm]
0, & x = \frac{1}{\pi n}, n \in \ZZZ.
\end{array}
\right.
$$
It is easy to verify that $\gfunc$ is $C^{\infty}$.
Then the function $\mrsfunc(x)=\int_{0}^{x}\gfunc(t)dt$ satisfies the conditions of Claim~\ref{clm:ex_not_in_J}.
\end{rem}
\proof
Suppose that $\mrsfunc = \dfval{\mrsfunc}{\Fld}$ for some smooth vector field $\Fld$.
Since $\pnt_{i}$ is a critical point, we have $d\mrsfunc(\pnt_{i})=0$ for all $i$, whence $\dfvalx{\mrsfunc}{\Fld}{\pnt_{i}}=0$.
On the other hand $\mrsfunc(\pnt_{i})\not=0$ by the assumption.
Hence the equality $\mrsfunc(\pnt_{i}) = \dfvalx{\mrsfunc}{\Fld}{\pnt_{i}}$ is impossible.
\endproof

\section{Proof of Theorem~\ref{th:Stab}}\label{sect:proofA}
\begin{lem}\label{lm:pStMP_in_Pcrf}
Suppose that $\mrsfunc\in\smone$ satisfies the condition \condBndCrVal.
Let also $\prMPP:\DiffM\times\DiffP\to\DiffP$ be the projection onto the second multiple.
Then 
$$
\begin{array}{lll}
\text{Case \ $\Psp=\RRR:$} & \qquad &  \prMPP(\StabfMR)\ \subset \ \DiffRfe \\ [1.5mm]
\text{Case \ $\Psp=\Circle:$} & \qquad &  \prMPP(\StabfMS)\ \subset \ \DiffSfE.
\end{array}
$$
\end{lem}
\proof
Let $(\difM,\difP)\in\StabfMS$, i.e., $\difP\circ\mrsfunc = \mrsfunc\circ\difM$.
Then $\difP =\prMPP(\difM,\difP)$.
We claim that $\difP$ preserves the set $\excf=\{1,\ldots,\nv\}$ of exceptional values of $\mrsfunc$.
Indeed, $\difM$ interchanges level-sets of $\mrsfunc$ in same manner as $\difR$ interchanges values of $\mrsfunc$.
Denote by $\singf$ the set of critical points of $\mrsfunc$.
Since $\difM$ is a diffeomorphism, it preserves the sets  $\mrsfunc^{-1}(\mrsfunc(\singf))$ and $\mrsfunc^{-1}(\mrsfunc(\partial\manif))$.
Therefore, $\difP$ preserves $\mrsfunc(\singf) \cup \mrsfunc(\partial\manif)=\excf$. 

This proves our lemma for the case $\Psp=\Circle$ as by the definition $\DiffSfE$ consists of preserving orientation diffeomorphisms that also preserves $\excf$.
If $\Psp=\RRR$, then $\difP$ also preserves the ordering of the {\em finite\/} set $\excf$ and therefore fixes it point-wise, i.e. $\difP\in\DiffRfe$.
\qed

Now Theorem~\ref{th:Stab} follows from Proposition~\ref{pr:proof_A} and~\ref{pr:proof_A_0} below.
\begin{prop}\label{pr:proof_A}
Suppose that $\nv\geq1$ and $\mrsfunc$ satisfies the conditions \condBndCrVal\ and \condJacId.
Then $\DiffPfe \subset \prMPP(\StabfMP)$ and $\prMPP$ admits a continuous section $$\sectMPP:\DiffPfe\to\StabfMP$$ being a homomorphism.
\end{prop}
\begin{prop}\label{pr:proof_A_0}
If $\nv=0$, then $\prMPP(\StabfMR)=\DiffRfe$ but a global section of $\prMPP$ exists if and only if $\mrsfunc:\manif\to\Circle$ is a trivial fibration.
\end{prop}
\subsection{Proof of Proposition~\ref{pr:proof_A}.}
The section $\sectMPP:\DiffPfe\to\StabfMP$ of the projection $\prMPP$ must have the following form
$\sectMPP(\difR) = (\difShift(\difR), \difR)$,
where $\difShift: \DiffPfe \to \DiffM$ is a continuous homomorphism such that
$\difR \circ \mrsfunc = \mrsfunc \circ \difShift(\difR).$
Thus we have to construct $\difShift$.

For $\cnv=1,\ldots,\nv$ let 

{\hangindent=1.5em  
$\bullet$~$\exlev{\cnv}=\mrsfunc^{-1}(\cnv)$ be the $\cnv$-th exceptional level of $\mrsfunc$,
}

{\hangindent=1.5em  
$\bullet$~$\lnbh{\cnv}=\mrsfunc^{-1}(\cnv,\cnv+1)$ be the part of $\manif$ between the levels $\exlev{\cnv}$ and $\exlev{\cnv+1}$,
}

{\hangindent=1.5em  
$\bullet$~$\elnbh{\cnv} = \mrsfunc^{-1}(\cnv-{\textstyle\frac{1}{3}}, \cnv+{\textstyle\frac{1}{3}})$ be a neighborhood of $\exlev{\cnv}$,
}

{\hangindent=1.5em  
$\bullet$~$\elnbhm{\cnv} = \mrsfunc^{-1}(\cnv-{\textstyle\frac{1}{3}}, \cnv)$ and 
 $\elnbhp{\cnv} = \mrsfunc^{-1}(\cnv, \cnv+{\textstyle\frac{1}{3}})$ respectively {\em lower\/} and {\em upper\/} part of $\elnbh{\cnv}$.
}

Thus the restrictions of $\mrsfunc$ to $\elnbh{\cnv}$ and to $\lnbh{\cnv}$ can be regarded as functions 
$$
\elnbh{\cnv}\to(\cnv-{\textstyle\frac{1}{3}}, \cnv+{\textstyle\frac{1}{3}}) 
\qquad \text{and} \qquad 
\lnbh{\cnv}\to(\cnv,\cnv+1). 
$$

\begin{lem}
In some neighborhood of $\exlev{\cnv}$ there is a vector field $\elFld{\cnv}$ such that \begin{equation}\label{equ:dfF_f_i}
\dfval{\mrsfunc}{\elFld{\cnv}}=\mrsfunc-\cnv.
\end{equation}
\end{lem}
\proof
Notice that for each point $\pnt\in\exlev{\cnv}$ there is a neighborhood $\elnbh{\pnt}\subset\elnbh{\cnv}$ and a vector field $\elFld{\pnt}$ defined on $\elnbh{\pnt}$ such that 
such that $\dfval{\mrsfunc}{\elFld{\pnt}}=\mrsfunc-\cnv$.
For critical points of $\mrsfunc$ this follows from the condition \condJacId\ and for the regular ones from (1) of Lemma~\ref{lm:class_JM}.
In the latter case we may put $\elFld{\cnv}=(x_1, 0,\ldots,0)$ in that coordinates for which $\mrsfunc(x_1,\ldots,x_{\dimM})=x_1+\cnv$.

Let $\{\elnbhref{\jj}\}$ be a finite refinement of the covering $\{\elnbh{\pnt}\}_{\pnt\in\exlev{\cnv}}$ of $\exlev{\cnv}$.
Then on each $\elnbhref{\jj}$ we have a vector field $\elFldref{\jj}$ such that 
$\dfval{\mrsfunc}{\elFldref{\jj}}=\mrsfunc-\cnv$.
Let $\mu_{\jj}:\elnbhref{\jj}\to[0,1]$ be a partition of unity subordinated to $\{\elnbhref{\jj}\}$, i.e.\!\! $\supp\mu_{\jj}\subset\elnbhref{\jj}$ and ${\textstyle\sum_{\jj}} \, \mu_{\jj} \equiv 1$.
Define the following vector field $\elFld{\cnv}$ in a neighborhood of $\exlev{\cnv}$ by 
$\elFld{\cnv} = \sum_{\jj} \mu_{\jj} \elFldref{\jj}.$
Then 
$$ 
\dfval{\mrsfunc}{\elFld{\cnv}}=
\dfval{\mrsfunc}{ \left({\textstyle\sum_{\jj}} \, \mu_{\jj}\, \elFldref{\jj}\right)}=
{\textstyle\sum_{\jj}}\, \mu_{\jj} \, \left( \dfval{\mrsfunc}{\elFldref{\jj}} \right)=
{\textstyle\sum_{\jj}}\, \mu_{\jj} \, (\mrsfunc-\cnv) = 
\mrsfunc-\cnv.
\qed
$$

Decreasing $\elnbh{\cnv}$ if necessary, we can assume that $\elFld{\cnv}$ is defined on some neighborhood of $\overline{\elnbh{\cnv}}$.
Then $\dfval{\mrsfunc}{\elFld{\cnv}}<0$ on $\elnbhm{\cnv}$ and 
$\dfval{\mrsfunc}{\elFld{\cnv}}>0$ on $\elnbhp{\cnv}$.
Let $\elflow{\cnv}:\elnbh{\cnv}\times(-\delta,\delta)\to\manif$ be a local flow generated by $\elFld{\cnv}$.

Notice, that $\mrsfunc$ has no critical points in $\lnbh{\cnv}$, whence there is a vector field $\lFld{\cnv}$ on $\lnbh{\cnv}$ such that $\dfval{\mrsfunc}{\lFld{\cnv}}>0$.
We can assume that $\lFld{\cnv}$ generates a {\em global\/} flow $\lflow{\cnv}:\lnbh{\cnv}\times\RRR\to\lnbh{\cnv}$.

Moreover, not violating condition~\eqref{equ:dfF_f_i} we can also suppose that 
$\elFld{\cnv}(\pnt)=\lFld{\cnv}(\pnt)$ for $\pnt\in\lnbh{\cnv}$,
and 
$\elFld{\cnv}(\pnt)=-\lFld{\cnv-1}(\pnt)$ for $\pnt\in\lnbh{\cnv-1}$
provided $\elFld{\cnv}(\pnt)$ is defined.

Then 
$$
\elflow{\cnv}(\pnt,t)=\lflow{\cnv}(\pnt,t) \ \text{for $\pnt\in\elnbhp{\cnv}$} 
\ \ \text{and} \ \
\elflow{\cnv}(\pnt,t)=\lflow{\cnv-1}(\pnt,-t) \ \text{for $\pnt\in\elnbhm{\cnv}$}.
$$

\begin{center}
\begin{figure}[ht]
\begin{tabular}{ccccc}
\includegraphics[height=5.5cm]{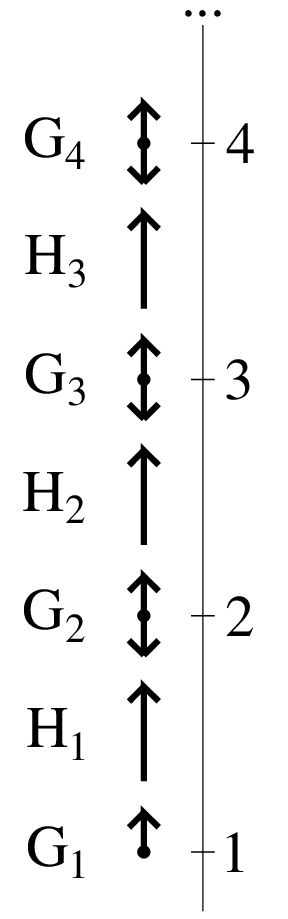} 
& \qquad \qquad &
\includegraphics[height=5.5cm]{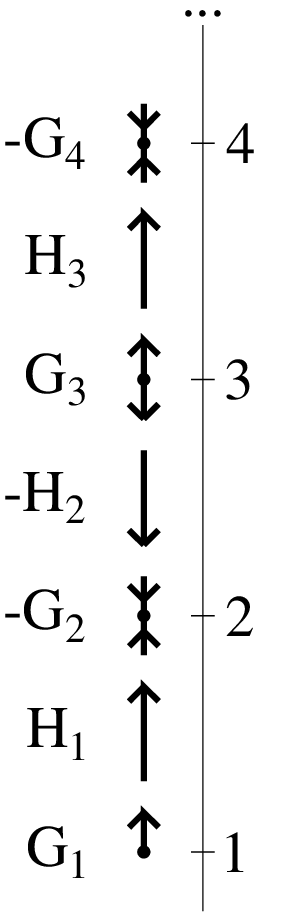}  
& \qquad \qquad &
\includegraphics[height=5.5cm]{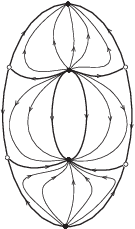} 
\\
a) & & b) & & c)
\end{tabular}
\caption{}\protect\label{fig:vf_R1}
\end{figure}
\end{center}

\begin{center}
\begin{figure}[ht]
\begin{tabular}{ccccc}
\includegraphics[width=3.8cm]{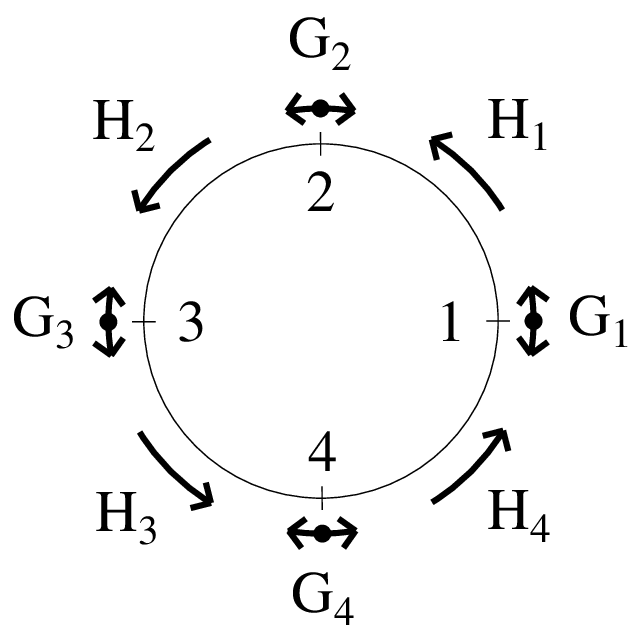} 
& \quad  &
\includegraphics[width=3.8cm]{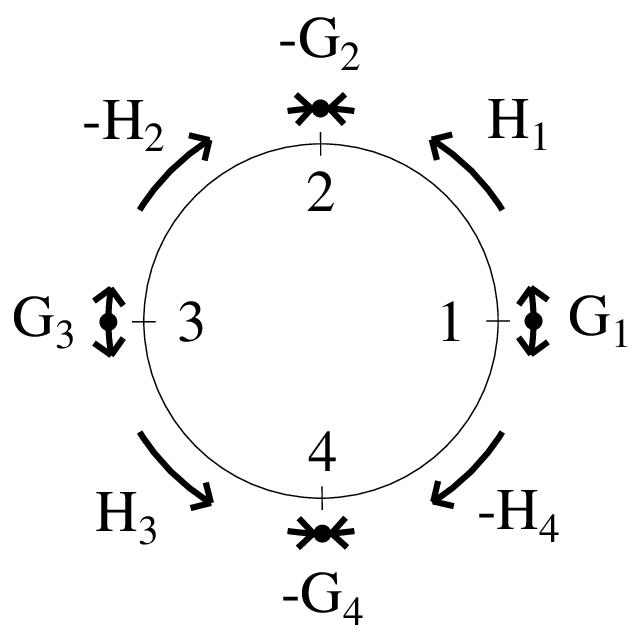} 
& \quad  &
\includegraphics[width=3.8cm]{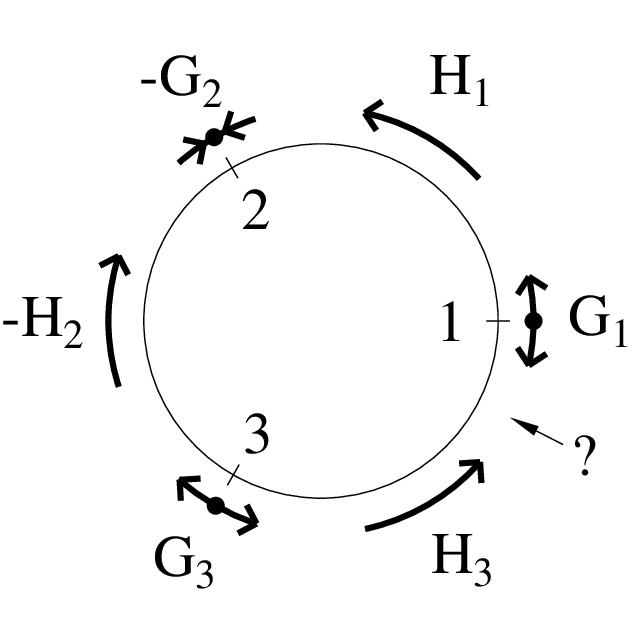} 
\\
a) & & b) & & c) 
\end{tabular}
\caption{}\protect\label{fig:vf_S1}
\end{figure}
\end{center}

\begin{rem}
Schematically vector fields $\elFld{\cnv}$ and $\lFld{\cnv}$ can be represented by the behavior of $\mrsfunc$ along their trajectories, see Figures~\ref{fig:vf_R1}a) and~\ref{fig:vf_S1}a).
Thus $\mrsfunc$ increases along $\lFld{\cnv}$ and along ``upper'' part of $\elFld{\cnv}$ and decreases along ``lower'' part of $\elFld{\cnv}$.
A bold point on the arrow means that $\elFld{\cnv}$ is tangent to the $\cnv$-th level-set of $\mrsfunc$.

Suppose that either $\Psp=\RRR$ or $\Psp=\Circle$ but $\nv$ is even. 
Then vector fields $\lFld{\cnv}$ and $\elFld{\cnv}$ give a global vector field $\rFld$ on $\manif$.
To construct $\rFld$ we should just change signs of all $\lFld{2\cnv}$ and $\elFld{2\cnv}$ (having even indices), see Figures~\ref{fig:vf_R1}b) and~\ref{fig:vf_S1}b).
In Figure~\ref{fig:vf_R1}c) such a vector field is shown for a height function on $2$-torus.
The critical points of $\mrsfunc$ are in bold. 
Together with white points they constitute the set of singular points of $\rFld$.
Existence of $\rFld$ would simplify the proof.
But we shall not use this approach since for the case $\Psp=\Circle$ and $\nv$ is odd such a vector field $\rFld$ does not exist, see Figure~\ref{fig:vf_S1}c).
\end{rem}


\begin{lem}\label{lm:ex_h}
For each $\difP\in\DiffPfe$ there exists a diffeomorphism $\difM\in\DiffM$ such that $\difP\circ\mrsfunc=\mrsfunc\circ\difM$, whence
$\difP\circ\mrsfunc\in\OrbfMP$ in fact belongs to $\OrbfM$.
\end{lem}
\proof
Let $\difP\in\DiffPfe$.
We shall find $C^{\infty}$ functions $\lfunc{\cnv}:\lnbh{\cnv}\to\RRR$ and $\elfunc{\cnv}:\elnbhref{\cnv}\to(-\delta,\delta)$, where $\elnbhref{\cnv}\subset\elnbh{\cnv}$ is a neighborhood of $\exlev{\cnv}$ such that the following mappings $\ldif{\cnv}:\lnbh{\cnv}\to\lnbh{\cnv}$ and $\eldif{\cnv}:\elnbhref{\cnv}\to\elnbh{\cnv}$ defined by: 
\begin{equation}\label{equ:def_ldif_eldif}
\ldif{\cnv}(\pnt)=\lflow{\cnv}(\pnt,\lfunc{\cnv}(\pnt)) 
\qquad \text{and} \qquad
\eldif{\cnv}(\pnt)=\elflow{\cnv}(\pnt,\elfunc{\cnv}(\pnt)), 
\end{equation}
are diffeomorphisms on their images, 
$\difP\circ\mrsfunc=\mrsfunc\circ\ldif{\cnv}$ on $\lnbh{\cnv}$, and 
$\difP\circ\mrsfunc=\mrsfunc\circ\eldif{\cnv}$ on $\elnbhref{\cnv}$.

Moreover, all these diffeomorphisms will coincide at common points of their ranges and therefore define a self-diffeomorphism $\difM$ of $\manif$ satisfying the statement of our lemma.

\begin{figure}[ht]
\includegraphics[height=3.5cm]{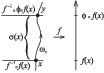}
\caption{}\protect\label{fig:def_sgm}
\end{figure}

{\bf Step 1. Definition of $\lfunc{\cnv}$ and $\ldif{\cnv}$}.
For each $x\in\lnbh{\cnv}$ let $\traj{\xpnt} \subset \lnbh{\cnv}$ be the trajectory of ${\xpnt}$ with respect to $\lflow{\cnv}$.
Notice that $\mrsfunc$ maps $\traj{\xpnt}$ diffeomorphically onto the open interval $(\cnv,\cnv+1)$ and $\traj{\xpnt}$ transversely intersects level-sets of $\mrsfunc$.
Let $y$ be a {\em unique\/} intersection point of $\traj{\xpnt}$ with the level-set
$\mrsfunc^{-1}(\difP\circ\mrsfunc(\xpnt))$ of $\mrsfunc$, see Figure~\ref{fig:def_sgm}.
Set $$\ldif{\cnv}(\xpnt)=y,$$
then $\mrsfunc\circ\ldif{\cnv}({\xpnt})=\difP\circ\mrsfunc(\xpnt)$.
Let also $\lfunc{\cnv}(\xpnt)$ be the time along $\traj{\xpnt}$ with respect to $\lflow{\cnv}$ from $\xpnt$ to $y$.
Then $y=\ldif{\cnv}(\xpnt)=\lflow{\cnv}(\xpnt,\lfunc{\cnv}(\xpnt))$.

Let us show that $\lfunc{\cnv}$ is smooth. This will imply a smoothness of $\ldif{\cnv}$.

Let $\pnt\in\lnbh{\cnv}$. 
Then $\lFld{\cnv}(\pnt)\not=0$, whence we can assume that in some local coordinates $(x_1,\ldots,x_{\dimM})$ at $\pnt$ we have $\lFld{\cnv}(x)=(1,0,\ldots,0)$.
For simplicity designate $\barx = (x_2,\ldots,x_{\dimM})$ and
$x=(x_1,\ldots,x_{\dimM})=(x_1,\barx)$.

It follows that $\lflow{\cnv}(x,t)=(x_1+t,\barx)$.
Suppose that a point $y$ belongs to the trajectory $\gamma_x$ of a point $x=(x_1,\barx)$.
Then $y=(y_1,\barx)$.
Moreover, the time between $x$ and $y$ on $\gamma_x$ is equal to $y_1-x_1$.

Since $\dfval{\mrsfunc}{\lFld{\cnv}} = \mrsfunc'_{x_1}\not=0$, it follows that there exists a unique smooth function $q(x)$
such that $x_1 = q(\mrsfunc(x), \barx)$.
Then from the definition of $\lfunc{\cnv}$ we get that
\begin{equation}\label{equ:clmdif_on_regpt}
\lfunc{\cnv}(x) =
q(\difP\circ\mrsfunc(x), \barx)-
q(\mrsfunc(x), \barx) =
q(\difP\circ\mrsfunc(x), \barx)- x_1,
\end{equation}
whence $\lfunc{\cnv}$ and $\ldif{\cnv}$ are smooth.

Let us verify that $\ldif{\cnv}$ is a diffeomorphism at each point $\pnt\in\lnbh{\cnv}$.
By~\cite[Theorem~19]{Maks:Shifts}, this is true if and only if the following relation holds true for all $\pnt\in\lnbh{\cnv}$:
\begin{equation}\label{d_clm_1}
\dfvalx{\lfunc{\cnv}}{\lFld{\cnv}}{\pnt}+1>0.
\end{equation}
Notice that 
$$
 \dfval{\lfunc{\cnv}}{\lFld{\cnv}} \ = \
 {\textstyle \frac{\partial \lfunc{\cnv}}{\partial x_1}} \ = \
q'_{x_1}(\difR\circ\mrsfunc(x), \barx) \cdot \difR'(\mrsfunc(x)) \cdot \mrsfunc'_{x_1}(x) \, - \, 1.
$$
Since $\difR$ is a preserving orientation diffeomorphism, we get $\difR'>0$.
Moreover, $\mrsfunc'_{x_1}$ and $q'_{x_1}$ have same signs.
Therefore
$\dfval{\lfunc{\cnv}}{\lFld{\cnv}} + 1 > 0$.

{\bf Step 2. Definition of $\elfunc{\cnv}$ and $\eldif{\cnv}$}.
Suppose that $\pnt\in\exlev{\cnv}$. 
Recall that we wish to define $\eldif{\cnv}$ by~\eqref{equ:def_ldif_eldif}.	
For simplicity assume that $\mrsfunc(\pnt)=0$.
Then $\difP(0)=0$ and by the Hadamard lemma $\difR(\vars)=\vars\cdot\bdifR(\vars)$ for some smooth function $\bdifR(\vars)$ such that $\bdifR(0)=\difR'(0)>0$.

Then by Lemma~\ref{lm:f_flow__f_exp}, see~\eqref{equ:main_relation}, the relation $\difP\circ\mrsfunc=\mrsfunc\circ\elflow{\cnv}(x,\elfunc{\cnv}(x))$ can be rewritten as follows:
$$
 \mrsfunc(x) \cdot \bdifR(\mrsfunc(x)) = \mrsfunc(x) \cdot e^{\eps \elfunc{\cnv}(x)},
$$
whence we set
\begin{equation}\label{equ:clmdif_on_crcomp}
\elfunc{\cnv}(x) = \eps \ln \bdifR(\mrsfunc(x)). 
\end{equation}
It follows that $\elfunc{\cnv}$ and $\eldif{\cnv}$ are smooth on some neighborhood $\elnbhref{\cnv}$ of $\exlev{\cnv}$.

Let us show that $\eldif{\cnv}$ is a diffemorphism, i.e. that $\dfvalx{\elfunc{\cnv}}{\elFld{\cnv}}{\pnt}>-1.$
Notice that
$$
\dfval{\elfunc{\cnv}}{\elFld{\cnv}} \ =  \
\dfval{(\eps \ln \bdifR(\mrsfunc))}{\elFld} \ = \
\frac{\eps \bdifR'_{t}(\mrsfunc)}{\bdifR(\mrsfunc)} \ \dfval{\mrsfunc}{\elFld} \ = \
\frac{\eps \bdifR'_{t}(\mrsfunc)}{\bdifR(\mrsfunc)} \ \mrsfunc.
$$
Since $\mrsfunc(\pnt)=0$, we get $\dfvalx{\elfunc{\cnv}}{\elFld{\cnv}}{\pnt}=0 > -1$.
Hence $\eldif{\cnv}$ is a diffeomorphism near $\pnt$.

{\bf Step 3. Coherency of $\ldif{\cnv}$ and $\eldif{\cnv}$.}
We have to show that $\eldif{\cnv}=\ldif{\cnv}$ on $\elnbhp{\cnv}$ and 
$\eldif{\cnv}=\ldif{\cnv-1}$ on $\elnbhm{\cnv}$. 

Let $\pnt\in\elnbhp{\cnv}$.
Then $\difP\circ\mrsfunc=\mrsfunc\circ\ldif{\cnv}(\pnt)=\mrsfunc\circ\eldif{\cnv}(\pnt)$. In other words:
$$
\mrsfunc\circ\lflow{\cnv}(\pnt,\lfunc{\cnv}(\pnt)) = 
\mrsfunc\circ\elflow{\cnv}(\pnt,\elfunc{\cnv}(\pnt)).
$$
Notice that $\mrsfunc$ is monotone along trajectories of $\lflow{\cnv}$ and those trajectories of $\elflow{\cnv}$ that does not belong to $\exlev{\cnv}$.
Since $\elflow{\cnv}(\pnt,t)=\lflow{\cnv}(\pnt,t)$ on $\elnbhp{\cnv}$, we get 
$\lfunc{\cnv}(\pnt)=\elfunc{\cnv}(\pnt)$ and 
$$
\eldif{\cnv}(\pnt)=
\elflow{\cnv}(\pnt,\elfunc{\cnv}(\pnt)) =
\lflow{\cnv}(\pnt,\lfunc{\cnv}(\pnt))=
\ldif{\cnv}(\pnt).
$$

Similarly, since $\elflow{\cnv}(\pnt,t)=\lflow{\cnv-1}(\pnt,-t)$ on $\elnbhm{\cnv}$, it follows that  
$\elfunc{\cnv}(\pnt)=-\lfunc{\cnv-1}(\pnt)$ and 
$\eldif{\cnv}(\pnt)=\elflow{\cnv}(\pnt,\elfunc{\cnv}(\pnt))=\lflow{\cnv-1}(\pnt,\elfunc{\cnv-1}(\pnt))=\ldif{\cnv-1}(\pnt)$.
Lemma~\ref{lm:ex_h} is proved.
\endproof

Thus the correspondence $\difP \mapsto \difM$ is a mapping
$\difShift:\DiffPfe\to\DiffM$ such that $(\difShift(\difP),\difP)\in\StabfMP$.
It remains to prove the following two statements.

\begin{claim}\label{clm:difShift_is homomorphism}
$\difShift$ is a homomorphism.
\end{claim}
\proof
Let $\difR_1,\difR_2\in \DiffPfe$, $\difR_0=\difR_1\circ\difR_2$,
$\difShift_{\kv} = \difShift(\difR_{\kv})$, $\kv=0,1,2$.
We have to show that $\hat\difShift=\difShift_1\circ\difShift_2\circ\difShift_0^{-1}$ coincides with $\id_{\manif}$.

First notice that $\hat\difShift\in\StabfM$. 
Indeed, since $\difR_{\kv}\circ\mrsfunc=\mrsfunc\circ\difShift_{\kv}$, we obtain that 
$$ 
\mrsfunc\circ\difShift_1\circ\difShift_2\circ\difShift_0^{-1} = \difR_1\circ\difR_2\circ\difR_0^{-1}\circ\mrsfunc =
\mrsfunc.
$$

Thus $\hat\difShift\in\StabfM$ and therefore it preserves each level-set of $\mrsfunc$.
Let $\pnt\in\lnbh{\cnv}$ and $\gamma_{\pnt}$ be the trajectory of $\pnt$ with respect to $\lflow{\cnv}$.
By the construction each $\difShift_{\kv}$ preserves $\gamma_{\pnt}$ for $\kv=0,1,2$, whence so does $\hat\difShift$.
Therefore $\hat\difShift$ also preserves the intersection 
$\pnt = \mrsfunc^{-1}\mrsfunc(\pnt)\cap \gamma_{\pnt}$, i.e. $\hat\difShift(\pnt)=\pnt$.
Thus $\hat\difShift$ is the identity on $\overline{\manif\setminus\exlevel}$.

If $\pnt\in\exlev{\cnv}$, then $\difShift_{\kv}(x)=\elflow{\cnv}(x,\elfunc{\cnv}^{\kv}(x))$ for some smooth function $\elfunc{\cnv}^{\kv}$, $\kv=0,1,2$.
But it follows from~\eqref{equ:clmdif_on_crcomp} that $\elfunc{\cnv}^{\kv}$ is constant on $\exlev{\cnv}$ as well as $\mrsfunc$.
This implies that $\hat\difShift$ is identity on $\exlev{\cnv}$ also.
\endproof

\begin{claim}
$\difShift:\DiffPfe \to \DiffM$ is continuous between strong $C^{\infty}$-topologies of these groups. 
\end{claim}
\proof
This follows from formulas~\eqref{equ:def_ldif_eldif}, \eqref{equ:clmdif_on_regpt}, and~\eqref{equ:clmdif_on_crcomp}.
We leave the details to the reader.
\endproof

\subsection{Proof of Proposition~\ref{pr:proof_A_0}.} 
Recall that $\manif$ is closed and $\mrsfunc:\manif\to\Circle$ is a locally trivial fibration.
Let us show that $\prMPP(\StabfMS)=\DiffSpl$.

Let $\lFld{}$ be a gradient vector field for $\mrsfunc$, then $\dfval{\mrsfunc}{\lFld{}}>0$ on all of $\manif$.
Let also $\lflow{}$ be the flow on $\manif$ generated by $\lFld{}$.
We claim that for every $\difR\in\DiffSpl$ there exists a smooth function $\lfunc{}:\manif\to\RRR$ such that the mapping  $\difM:\manif\to\manif$ defined by $\difM(x)=\lflow{}(x,\lfunc{}(x))$ is a diffeomorphism and $\difP\circ\mrsfunc=\mrsfunc\circ\difM$.

To construct $\lfunc{}$, notice that there exists an isotopy $\difP_t:\Circle\to\Circle$ such that $\difR_0=\id_{\Circle}$ and $\difR_1=\difR$.
Let $x\in\manif$, and $\gamma:[0,1]\to\Circle$ be the path of the point $\mrsfunc(x)$ under $\difR_t$, i.e. $\gamma(t)=\difP_t(\mrsfunc(x))$.
Since $\mrsfunc$ is a locally trivial fibration, we see that there exists a {\em unique\/} lifting $\omega:[0,1]\to\manif$ of $\gamma$ to $\manif$ that starts at $x$ and lies in the trajectory of $x$ with respect to $\lflow{}$. 
Thus $\omega(0)=x$, $\mrsfunc(\omega(x))=\gamma(t)=\difP_t(\mrsfunc(x))$, and $\omega[0,1]$ is a part of $\lflow{}$-trajectory of $x$.

Let $\lfunc{}(x)$ be the time along $\omega$ with respect to $\lflow{}$ and $\difM(x)=\lflow{}(x,\lfunc{}(x))$.
Then similarly to Step 2 of Lemma~\ref{lm:ex_h} it can be shown that $\lfunc{}$ and $\lflow{}$ are smooth.

Notice that the definition of $\lfunc{}$ and $\difM$ depends on the isotopy $\difP_t$.
Since $\DiffSpl$ is not contractible, it follows that in general we can not choose $\difM$ to continuously depend on $\difR$.

Suppose $\mrsfunc$ is a trivial fibration, i.e. $\manif = \nmanif\times\Circle$ and $\mrsfunc:\nmanif\times\Circle\to\Circle$ is given by $\mrsfunc(x,\vars)=\vars$.
Then the homomorphism $\difShift:\DiffSpl\to\DiffM$ can be defined by 
$\difShift(\difP)(x,\vars) = (x,\difP(\vars))$ and we also have that 
\ $\difP\circ\mrsfunc(x,\vars) = \difP(\vars) = \mrsfunc\circ\difShift(\difP)(x,\vars).$ 

Conversely, suppose that there exists a continuous mapping $\difShift:\DiffSpl\to\DiffM$ such that $\difP\circ\mrsfunc=\mrsfunc\circ\difShift(\difP)$.
We can always assume that $\difShift(\id_{\Circle})=\id_{\manif}$.
Let $\difP_t:\Circle\to\Circle$ be the ``rotation'' isotopy of $\Circle$:
$\difP_t(\vars)=\vars+t \mod 1$.
This isotopy yields an isotopy $\difM_t=\difShift(\difP_t)$ of $\manif$
($\difM_t$ is \`a priori just continuous, but this is enough for the proof of triviality of $\mrsfunc$).
Since $\difP_0=\difP_1=\id_{\Psp}$, we also have that $\difM_0=\difShift(\difP_0)=\difShift(\difP_1)=\difM_1$.

Denote $\nmanif=\mrsfunc^{-1}(0)$. 
Then $\difM_t$ yields a continuous map $\qfunc:\nmanif\times\Circle\to\manif$ defined by $\qfunc(x,\vars)=\difM_{\vars}(x)$.
It is easy to verify that $\qfunc$ is a homeomorphism.
Moreover, the composition 
\ $\nmanif\times\Circle \xrightarrow{\ \qfunc\ }  \manif \xrightarrow{\ \mrsfunc\ } \Circle$ \
is given by
$$
\mrsfunc\circ\qfunc(x,\vars)=
\mrsfunc\circ\difM_{\vars}(x)=
\difP_{\vars}\circ\mrsfunc(x)=
\mrsfunc(x)+\vars = 0+\vars=\vars,
$$
i.e. it coincides with the projection $\nmanif\times\Circle\to\Circle$.
Hence $\mrsfunc$ is a trivial fibration.
\qed\section{Factor spaces $\DiffRf/\DiffRfe$ and $\DiffSpl/\DiffSfe$.}\label{sect:DiffPfe}
Let us endow the spaces $\DiffR$ and $\DiffSpl$ with strong $C^{r}$ Whitney topologies for some $r=0,1,\ldots,\infty$.
These topologies yields certain topologies on $\DiffRf$, $\DiffRfe$, $\DiffSfe$, and on the spaces of \emph{left} adjacent classes $\DiffRf/\DiffRfe$ and $\DiffSpl/\DiffSfe$.

\begin{lem}\label{lm:DiffP_is_contr}
The groups $\DiffRf$, $\DiffRfe$, and $\DiffSfe$ for $\nv\geq1$ are contractible.
\end{lem}
\proof
For each of these groups its contracting $H$ can be defined by the standard formula:
$H(\difR,t)(x) = (1-t) x + t \difR(x).$
\endproof

\begin{thm}\label{th:struct_DP_DPe}
There are homeomorphisms
$$
\DiffRf/\DiffRfe \approx \RRR^{\nv-2} \qquad \text{and} \qquad
\DiffSpl/\DiffSfe \approx \Circle\times\RRR^{\nv-1}.
$$
\end{thm}
\proof[Proof. Case $\Psp=\RRR$.]
Let $\padj:\DiffRf \to \DiffRf/\DiffRfe$ be the natural projection.
Consider the following subset of $\RRR^{\nv-2}$:
$$
\dntwo= \{ (x_{2},\ldots,x_{\nv-1}) \in \RRR^{\nv-2} \ | \
1<x_2 < x_3 < \cdots < x_{\nv-2}< n \}.
$$
Evidently, $\dntwo$ is open and convex subset of $\RRR^{\nv-2}$, whence it is diffeomorphic with $\RRR^{\nv-2}$.
We shall show that $\DiffRf/\DiffRfe$ is canonically homeomorphic with $\dntwo$.

Consider the \emph{evaluation map} $\pcr:\DiffRf\to\dntwo$ defined by the formula $\pcr(\difR)=\bigl(\difR(2), \ldots, \difR(\nv-1)\bigr)$ for $\difR\in\DiffRf$.
This mapping evidently factors to a unique bijection $\adjd:\DiffRf/\DiffRfe\to\dntwo$ such that
$$
\begin{CD}
\pcr=\adjd\circ\padj:\DiffRf @>{\padj}>>\DiffRf/\DiffRfe @>{\adjd}>> \dntwo.
\end{CD}
$$

Our aim is to prove that $\adjd$ is a homeomorphism for all $C^{r}$-topologies.
It is easy to see that $\pcr$ is continuous in $C^{0}$-topology ($C^{0}$-continuous).
Hence it is so in all $C^{r}$-topologies for $r=1,\ldots,\infty$.
Then from the definition of the factor-topology on $\DiffRf/\DiffRfe$, we obtain that $\adjd$ is also $C^{r}$-continuous for $r=0,\ldots,\infty$.

We will now show that $\pcr$ admits a section $\sct:\dntwo\to\DiffRf$ which is $C^{r}$-continuous for every $r=0,\ldots,\infty$.
It follows that so is $\adjd^{-1}=\pcr\circ\sct$, whence $\adjd$ is a $C^{r}$-homeomorphism.

Consider another subset of $\RRR^{\nv}$:
$$\dn= \{ (x_{1},\ldots,x_{\nv}) \in \RRR^{\nv} \ | \
0< x_1< x_2  < \cdots < x_{\nv}< \nv+1 \}.
$$
Then $\dntwo$ can be identified with the subset
$\{1\} \times \dntwo \times \{\nv\}$ of $\dn$.

\begin{lem}\label{clm:constr_section_not_onto}
There exists a smooth function $\mnfunc:\dn\times\RRR \to \RRR$
with the following properties:
\begin{enumerate}
\em\item
$\mnfunc'_{t}(x_1,\ldots,x_{\nv}, t) >0$ for all $(x_1,\ldots,x_{\nv}, t) \in\dn\times\RRR$;
\item
$\mnfunc(x_1,\ldots,x_{\nv}, t) = t$ for $t\leq 0$ and $t\geq \nv+1$;
\item
$\mnfunc\left( x_1,\ldots,x_{\nv}, k \right) = x_k$ for $k=1,2,\ldots,\nv$.
\item
$\mnfunc\left( 1,2,\ldots,\nv, t \right) = t$ for $t\in\RRR$.
\end{enumerate}
\end{lem}
It follows from this lemma that $\mnfunc$ yields a section  
$\sct:\dntwo\to\DiffRf$ defined by the following formula:
\begin{equation}\label{equ:sect}
\sct(x_2,\ldots,x_{\nv-1})(t) = \mnfunc(1,x_2,\ldots,x_{\nv-1},\nv, t).
\end{equation}

Notice that by (2) all diffeomorphisms $\sct(x_2,\ldots,x_{\nv-1})$ are fixed outside $[0,\nv+1]$.
This implies that $\mnfunc$ is uniformly continuous, whence $\sct$ is continuous in all $C^{r}$-topologies.

Indeed, let $x=(x_2,\ldots,x_{\nv-1}) \in \dntwo$, 
$\eps:\manif\to(0,\infty)$ be a strictly positive continuous function,
and $B_{\eps}$ be a base neighborhood of $\difR$ in some $C^{r}$-topology consisting of $\psi\in\DiffRf$ such that
$\mathop\sum\limits_{i=0}^{r}\|D^{k}\sct(x)(z)-D^{k}\psi(z)\| < \eps(z)$ for all $z\in\manif$.

By uniform continuity of $\mnfunc$ there is a neighborhood $V$ of $x$ in $\dn$ such that $\mathop\sum\limits_{i=0}^{r}|D^{k}\sct(x)(z)-D^{k}\sct(y)(z)|<\eps(z)$
for all $y\in V$ and $z\in\manif$.
Thus $\sct(V\cap\dntwo) \subset B_{\eps}$, whence $\sct$ is continuous in $C^{r}$-topology for all $r<\infty$.
Therefore it is also continuous in $C^{\infty}$-topology.
This proves case $\RRR$ of our theorem modulo Lemma~\ref{clm:constr_section_not_onto}.
\endproof

\proof[Proof of Lemma~\ref{clm:constr_section_not_onto}]
Suppose that $\mnfunc$ is defined.
Set $\dmnfunc=\mnfunc'_{t}$.
Then 
\begin{enumerate}
\item[(1$'$)]
$\dmnfunc>0$;
\item[(2$'$)]
$\dmnfunc(x_1,\ldots,x_{\nv}, t)=1$, for $t\not\in(0,\nv+1)$;
\item[(3$'$)]
$\int\limits_{k}^{k+1}\dmnfunc(x_1,\ldots,x_{\nv}, t) dt = x_{k+1} - x_{k}$;
\item[(4$'$)]
$\dmnfunc(1,2,\ldots,\nv, t) \equiv 1$,
\end{enumerate}
and $\mnfunc(x,t) = \int\limits_{0}^{t}\dmnfunc(x,s)ds$.

Thus it suffices to construct a function $\dmnfunc$ satisfying the conditions (1$'$)-(4$'$).
Then (1)-(4) will be satisfied.

For every pair $a < b \in\RRR$ we will now define a smooth function $\qfunc_{a,b}(t,s)\geq 0$ such that
\begin{enumerate}
\item[(2$''$)]
$\qfunc_{a,b}(t,s)=0$ for $t\not\in(a,b)$ and $s\in\RRR$, and
\item[(3$''$)]
$\int\limits_{a}^{b}\qfunc_{a,b}(t,s)dt = s-(b-a)$.
\end{enumerate}

It follows from (3$''$) and the condition $\qfunc_{a,b}(t,s)\geq 0$ that
\begin{enumerate}
\item[(4$''$)]
$\qfunc_{a,b}(t,b-a)=0$ for $t\in\RRR$.
\end{enumerate}

Then $\dmnfunc$ can be defined as follows:
\begin{multline*}
\dmnfunc(x_1,\ldots,x_n, t) =
1 + \sum_{k=0}^{\nv}\qfunc_{k,k+1}(t,x_{k+1}-x_{k}) = \\
= 1 + \qfunc_{0,1}(t,x_1) \ + \ \qfunc_{1,2}(t,x_2-x_1) \ +
 \ \cdots \ + \ \qfunc_{\nv,\nv+1}(t,x_{\nv+1}-x_{\nv}).
\end{multline*}

Indeed, $\dmnfunc>0$.
Further, (2$'$) and (4$'$) follows from (2$''$) and (4$''$) respectively.
Finally, let us verify (3$'$). If $k=0,\ldots,\nv$, then we have

\begin{multline*}
\int\limits_{k}^{k+1}\dmnfunc(x_1,\ldots,x_{\nv}, t) dt
\stackrel{(2'')}{=}
\int\limits_{k}^{k+1} \left( 1 + \qfunc_{k,k+1}(t,x_{k+1}-x_{k}) \right) dt
\stackrel{(3'')}{=} \\ =
(k+1 - k) + (x_{k+1}-x_{k} - 1) = x_{k+1}-x_{k}.
\end{multline*}

Thus it remains to construct $\qfunc_{a,b}\geq0$ satisfying (2$''$) and (3$''$).
Consider the following $C^{\infty}$-functions
\begin{gather*}
\afunc(t) = \left\{
\begin{array}{cl}
\frac{1}{b-a} \cdot e^{- \frac{1}{(t-a)(b-t)}}, & t \in (a,b) \\
0, & t \not\in (a,b),
\end{array}
\right.
\\
 \bfunc(t,c) = e^{c \cdot \afunc(t)}-1,
\quad \qquad \text{and} \qquad \quad
 \cfunc(c) = \int\limits_{a}^{b} \bfunc(t,c) dt.
\end{gather*}
Then $\cfunc(c)$ is smooth,
$\lim\limits_{c\to-\infty}\cfunc(c)=-(b-a)$,
$\lim\limits_{c\to+\infty}\cfunc(c)=+\infty$, $\cfunc(0)=0$, and
$$
\cfunc'_{c}(c) =
\int\limits_{a}^{b} \bfunc'_{c}(t,c) dt =
\int\limits_{a}^{b} \afunc(t) e^{c \cdot \afunc(t)} dt > 0.
$$

Thus $\cfunc$ is a diffeomorphism of $\RRR$ onto $\bigl(-(b-a),+\infty\bigr)$.

Let $\cfuncinv:\bigl(-(b-a),\infty\bigr) \to \RRR$ be a diffeomorphism inverse to $\cfunc$.
Then the following function satisfies (2$''$) and (3$''$):
$$ \qfunc_{a,b}(t,s) = \bfunc(t,\cfuncinv(s-(b-a))).$$
Indeed, (2$''$) is obvious.
Let us verify (3$''$):
\begin{multline*}
\int\limits_{a}^{b} \qfunc_{a,b}(t,s) dt \ = \
\int\limits_{a}^{b} \bfunc(t,\cfuncinv(s-(b-a))) dt \ = \ 
\cfunc(\cfuncinv(s-(b-a))) \ = \ s-(b-a). \qed
\end{multline*}

\proof[Proof. Case $\Psp=\Circle$.]
Recall that the $\nv$-th configuration space of $\Circle$ is the following subset 
\[
\ConfSp=\{ (x_1,\ldots,x_{\nv}) \in\Tn \ | \ x_i\not=x_j \ 
\text{for}\ i\not=j \}
\]
of $n$-dimensional torus $\Tn=\Circle\times\cdots\times\Circle$.

Let $\impcr$ be the connected component of $\ConfSp$ containing the point 
$(1,\ldots,\nv)$.
Denote also 
\[
\dnone= \{ (x_{2},\ldots,x_{\nv}) \in \RRR^{\nv-1} \ | \
0 < x_2 < x_3  < \cdots < x_{\nv} < \nv \}.
\]
Then $\dnone$ is an open and convex subset of $\RRR^{\nv-1}$, whence it is diffeomorphic with $\RRR^{\nv-1}$.

\begin{lem}\label{lm:Fn_S1_symplex}
$\impcr$ is {\em diffeomorphic\/} with $\dnone\times\Circle$.
\end{lem}
\proof
Recall that we regard $\Circle$ as $\RRR/\nv\ZZZ$.
Consider the following mapping $\hconf:\ConfSp\to\dnone\times\Circle$ defined by 
\[
\hconf(x_{1}, \ldots, x_{\nv-1}, x_{\nv}) =
\bigl(
 x_{1}-x_{\nv}, \ \ldots, \ x_{\nv-1}-x_{\nv}, \ [x_{\nv}] 
\bigr),
\]
where the differences are taken modulo $\nv$, and $[x]$ means $x\!\!\mod\nv$.
Since $x_{\acnv}\not=x_{\acnv'}$ for $\acnv\not=\acnv'$, it follows that $\hconf$ is well defined.
Notice also that $\hconf(1,\ldots,\nv-1,\nv)=(1,\ldots,\nv-1,[\nv])$.
Evidently, $\hconf$ is smooth and admits a smooth right inverse $s:\dnone\times\Circle\to\ConfSp$ defined by 
$s(t_1,\ldots,t_{\nv-1},x)=(t_1+x, \ldots, t_{\nv-1}+x, x).$
Since $\dnone\times\Circle$ is connected, it follows that $\hconf$ yields a diffeomorphism of the connected component of $\ConfSp$ containing the point $(1,\ldots,\nv)$, i.e.\!\! $\impcr$, onto $\dnone\times\Circle$.
\endproof

Let $\pcr:\DiffSf\to\ConfSp$ be the \emph{evaluation map\/} defined by
\begin{equation}\label{equ:eval_map}
\pcr(\difR)=\left( 
\difR(1), \ldots, \difR(\nv) \right),
\quad \difR\in\DiffSf.
\end{equation}
Since $\DiffSf$ is connected, it follows that the image $\pcr(\DiffSf)$ coincides with the connected component $\impcr$ of $\ConfSp$ containing the point $\pcr(\id_{\Circle})=(1,\ldots,\nv)$.

Again $\pcr$ is constant on the adjacent classes of $\DiffSf$ by $\DiffSfe$ and yields a bijective continuous mapping $\tpcr:\DiffSf/\DiffSfe\to\impcr$ such that
$$
\begin{CD}
\pcr=\tpcr\circ\padj:\DiffSf @>{\padj}>>\DiffSf/\DiffSfe @>{\tpcr}>> \impcr,
\end{CD}
$$
where $\padj$ is a factor-mapping.
Thus in order to show that $\tpcr$ is a homeomorphism it suffices to prove that
$\hconf\circ\tpcr:\DiffSf/\DiffSfe\to\dnone\times\Circle$ is a homeomorphism in all Whitney topologies.

Consider the composition:
$$
 \Pcr =\hconf\circ\tpcr\circ\padj:
 \DiffSf \ \xrightarrow{\ \padj \ } \
 \DiffSf/\DiffSfe \ \xrightarrow{\ \tpcr\ } \
 \impcr \ \xrightarrow{\ \hconf\ } \
 \dnone\times\Circle
$$
As in the previous case it suffices to find a continuous section $\sct:\dnone\times\Circle\to\DiffSf$ of $\Pcr$.

Let $\mnfunc:\dnone\times\RRR\to\RRR$ be the function constructed in Lemma~\ref{clm:constr_section_not_onto} but for $\nv-1$. 
Then $\mnfunc$ maps the set $\dnone\times[0,\nv]$ onto $[0,\nv]$ so that
$\mnfunc(x,0)=0$ and $\mnfunc(x,\nv)=\nv$ for all $x\in\dnone$.
Hence $\mnfunc$ yields a \emph{continuous\/} mapping
$\tmnfunc:\dnone\times\Circle\to\Circle$ defined by factorization of $[0,\nv]$ onto $\Circle$ by gluing the ends $0$ and $\nv$.

Moreover,
$\mnfunc'_{t}(x,0)=\mnfunc'_{t}(x,\nv)=1$ and
$\mnfunc^{(s)}_{t}(x,0)=\mnfunc^{(s)}_{t}(x,\nv)=0$ for $s\geq 2$ and $x\in\dnone$.
Hence, $\tmnfunc$ is in fact $C^{\infty}$.
Then it is easy to verify that the following mapping $\sct:\dnone\times\Circle\to\DiffSf$ defined by:
\begin{equation}\label{equ:sectS1}
\sct(x_1,\ldots,x_{\nv-1},t)(\tau) \ = \ \tmnfunc(x_1,\ldots,x_{\nv-1},\tau)+t \mod \nv
\end{equation}
is a continuous section of $\Pcr$.
Theorem~\ref{th:struct_DP_DPe} is completed. 
\endproof

\subsection{Cyclic actions on $\ConfSp$.}
Consider the action of the group $\ZZZ_{\nv}$ on $\ConfSp$ by cyclic shift 
$\genzn:\ConfSp\to\ConfSp$ of coordinates:
$$
\genzn(x_{1},\ldots,x_{\nv})=(x_{2},\ldots,x_{\nv},x_{1}).
$$
Evidently, this action is free and $\genzn(\impcr)=\impcr$.
Moreover, $\genzn$ preserves orientation of $\impcr$ iff and only if $\nv$ is odd.

Suppose that $\nv=\kcnv\dcnv$ and let $\Zk$ be a cyclic subgroup of order $\kcnv$ of $\ZZZ_{\nv}$ generated by $\genzn^{\dcnv}$.

\begin{lem}\label{lm:FnZk_S1_symplex}
Suppose that $\nv$ is even and $\dcnv=\nv/\kcnv$ is odd, then $\impcr/\Zk$ is diffeomorphic with a ``twisted'' product $\dnone\,\widetilde{\times}\,\Circle$.
Otherwise, we have $\impcr/\Zk\approx\dnone\times\Circle$.
\end{lem}
\proof
The proof is direct and based on the remark that 
$\genzn^{\dcnv}$ reverses orientation if and only if $\nv$ is even and $\dcnv$ is odd.
\endproof
\section{Exceptional values}\label{sect:cont_crvmap}
Suppose that $\mrsfunc\in\smone$ satisfies the condition \condBndCrVal.
To each $\gfunc\in\OrbfMP$ we will now correspond the set of its exceptional values and give the condition on $\mrsfunc$ when this correspondence is continuous.

If $\Psp=\RRR$, then for every $\gfunc\in\OrbfMR$ the ordered (by increasing) set of its exceptional values represents a point in 
$$\dntwo\equiv\{1\}\times\dntwo\times\{\nv\}\subset\dn,$$
 since the minimal and maximal values $1$ and $\nv$ are fixed with respect to $\DiffMRf$.
Hence we have a well-defined mapping $$\crvmap:\OrbfMR\to\dntwo.$$

Suppose that $\Psp=\Circle$.
Then the exceptional values of $\mrsfunc$ are ordered only cyclically.
Therefore the set of exceptional values of $\gfunc\in\OrbfMS$ is a point $[\gfunc]$ in the factor space $\ConfSp/\ZZZ_{\nv}$.
Moreover, since we act by the {\em connected\/} group $\DiffSf$, we see that $[\gfunc]$ belongs to the connected component of $\ConfSp/\ZZZ_{\nv}$ containing the class of a set ${1,\dots,\nv}$ of exceptional values of $\mrsfunc$.
This connected component is $\impcr/\ZZZ_{\nv}$, whence we get the following mapping: $$\crvmap:\OrbfMS\to\impcr/\ZZZ_{\nv}.$$
By Lemma~\ref{lm:FnZk_S1_symplex} $\impcr/\ZZZ_{\nv}$ is $\Circle\times\dnone$ if $\nv$ is odd and $\Circle\widetilde{\times}\dnone$ for even $\nv$.

We will now give a sufficient conditions for continuity of $\crvmap$ and show that without this condition $\crvmap$ can loose continuity even in $C^{\infty}$-topology of $\OrbfMP$.

\begin{lem}\label{lm:cond_crlev_cont}
Suppose that every {\em critical\/} level-set of $\mrsfunc$
includes either a connected component of $\partial\manif$ or an essential critical point of $\mrsfunc$.
Then the mapping $\crvmap$ is continuous in $C^{\infty}$-topology of $\OrbfMP$.
\end{lem}
\proof
It suffices to prove continuity of $\crvmap$ at $\mrsfunc$.
Choose some $\eps\in(0,1/3)$ and let
$\cnbh=\mathop\cup\limits_{i=1}^{\nv}(i-\eps, i+\eps)$
be a neighborhood of the set of exceptional
values of $\mrsfunc$.
We have to find a $C^{\infty}$-neighborhood
$\Nbh$ of $\mrsfunc$ in $\smone$ such that for $\crvmap(\Nbh\cap\OrbfMP) \subset \cnbh$.

For every $i=1,\ldots,\nv$ let $\crset_{i} = \singf \cap \mrsfunc^{-1}(i)$
be the set of critical points belonging to
the $i$-th exceptional level-set of $\mrsfunc$.
Evidently, $\crset_{i}\cap\crset_{j}=\emptyset$ for $i\not=j$.
Also $\crset_{i}=\emptyset$, iff $i$ is a boundary but not a critical value of $\mrsfunc$.

We will assume that $\manif$ is given with some Riemannian metric.
Then for every $\gfunc\in\smone$ the norm $\|d\gfunc(x)\|$ of the differential of $\gfunc$ at every $x\in\manif$ is well defined.
Notice that $\|d\mrsfunc\|=0$ on $\crset_{i}$.
Hence there is a $\delta>0$ and for every $i=1,\ldots,\nv$ a compact neighborhood $\anbh_{i}$ of $\crset_{i}$ such that

\begin{itemize}
\item
$\anbh_{i}=\emptyset$ provided $\crset_{i}=\emptyset$;
\item
$\|d\mrsfunc\|>\delta$ on $\overline{\manif\setminus\anbh}$,
where $\anbh=\mathop\cup\limits_{i=1}^{n}\anbh_{i}$, and

\item
$\anbh_{i}\cap\anbh_{j}=\emptyset$ for $i\not=j$.
\end{itemize}

Let $\Lambda$ be the set of critical values of $\mrsfunc$ that are not boundary ones.
If $i\in\Lambda$, then by the assumptions there is an essential critical point
$\pnt_{i}\in\crset_{i} \subset\anbh_{i}$.
Hence for a neighborhood $\anbh_{i}$ we can find a neighborhood $\Nbh_{i}$ of
$\mrsfunc$ in $\smone$ such that every $\gfunc\in\Nbh_{i}$ has a critical point in $\anbh_{i}$.

Let also $\Nbh_{0}$ be a $C^{1}$-neighborhood of $\mrsfunc$ in $\smone$ consisting of
smooth functions $\gfunc$ such that
\begin{enumerate}
\item[{\rm(i)}]
$\|d\gfunc\|>\delta$ on $\overline{\manif\setminus\anbh}$, and

\item[{\rm(ii)}]
$|\mrsfunc-\gfunc|<\eps$.
\end{enumerate}

Denote $\Nbh=\Nbh_{0}\cap\bigl(\,\mathop\cap\limits_{i\in\Lambda}\Nbh_{i}\,\bigr)$.
We claim that $\crvmap(\Nbh\cap\OrbfMP) \subset \cnbh$.

Indeed, suppose that $\gfunc \in \Nbh\cap\OrbfMP$.
We have to show that the $i$-th exceptional value of $\gfunc$ differs from $i$ less than $\eps$.
Since $\gfunc$ has precisely $\nv$ exceptional values, it suffices to prove that every interval of the form $(i-\eps,i+\eps)$ for $i=1,\ldots,\nv$ contains an exceptional value of $\gfunc$.

It follows from {\rm(i)} that $\gfunc$ may have critical points only in $\anbh$.
If $i\in\Lambda$, then $\gfunc$ has a critical point in $\anbh_{i}$ by the construction of $\Nbh_{i}$.

Moreover, by {\rm(ii)} the boundary values of $\gfunc$ on the connected components of $\partial\manif$ differ from the corresponding values of $\mrsfunc$ less than $\eps$.
Thus every interval $(i-\eps,i+\eps)$ for $i=1,\ldots,\nv$ contains an exceptional value of $\gfunc$.
Lemma is proved.
\endproof

\subsection{Example of a function for which $\crvmap$ is not continuous.}
Let $\mrsfunc:\RRR\to\RRR$ be a smooth non-decreasing function having precisely two {\em flat\/} critical points $a<b$, i.e. $\mrsfunc^{(r)}(a)=\mrsfunc^{(r)}(b)=0$ for $r\geq 1$.

\begin{figure}[ht]
\includegraphics[height=4cm]{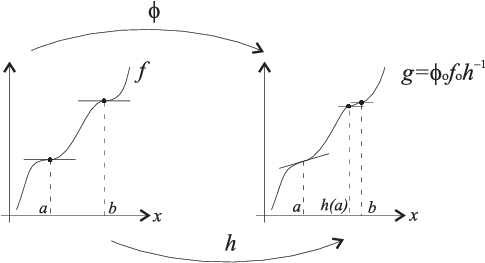}
\caption{}\protect\label{fig:k_is_not_cont}
\end{figure}

Let $\nbh_{a}$ and $\nbh_{b}$ be two disjoint $\eps$-neighborhoods of $a$ and $b$ for some $\eps>0$, $\nbh$ be a neighborhood of $[a,b]$,
and $0 < \delta < \frac{1}{3}( \mrsfunc(b)-\mrsfunc(a) )$.

Let $\aNbh$ be a neighborhood of $\mrsfunc$ in $C^{\infty}(\RRR)$ consisting of functions $\gfunc$ such that \ 
$\|g-f\|_{\overline{\nbh}}^{r} \ = \ \sup\limits_{x\in\overline{\nbh}}\
 \mathop\sum\limits_{i=0}^{r} |\gfunc^{(r)}(x) - \mrsfunc^{(r)}(x)| \ < \ \delta$.

It is easy to see that $\aNbh$ contains a function $\gfunc=\difR\circ\mrsfunc\circ\difM^{-1}\in\OrbfMR$, where $\difR$ and $\difM$ are diffeomorphisms of $\RRR$ such that $\difM(a)\in\nbh_{b}$, $\difM(b)=b$, $\gfunc\circ\difM(a)=\difR\circ\mrsfunc(a)\in
(\mrsfunc(b)-\delta,\mrsfunc(b)+\delta)$, and
$\gfunc\circ\difM(b)=\difR\circ\mrsfunc(b)=\mrsfunc(b)$.

The construction of $\gfunc$ is illustrated on the right diagram of Figure~\ref{fig:k_is_not_cont}.
The idea is that a critical point $a$ can be eliminated by arbitrary small perturbation of $\mrsfunc$ in any $C^{r}$-topology.
On the other hand, we can ``create'' a critical point equivalent to $a$ in arbitrary small neighborhood of $b$.
All this can be produced by conjugating $\mrsfunc$ so that 
$\|g-f\|_{\overline{\nbh}}^{r} < \delta$.

Thus both critical values $\gfunc\circ\difM(a)$ and $\gfunc\circ\difM(b)$ of $\gfunc$ belong to $\nbh_{b}$, i.e. the distance $\RRR^2$ between points $(\gfunc\circ\difM(a),\gfunc\circ\difM(b))$ and
$(\mrsfunc(a),\mrsfunc(b))$ is greater that $\delta$.

Since this holds for arbitrary small $\delta$, we obtain that $\crvmap$ is not continuous in $C^{r}$-topology for arbitrary $r\geq0$, whence it is not continuous in $C^{\infty}$-topology.
 
\section{Proof of Theorem~\ref{th:Orbit}}\label{sect:proofB}
Suppose that $\mrsfunc\in\smone$ satisfies the conditions \condJacId\ and \condBndCrVal, and 
every critical level-set of $\mrsfunc$ either has an essential critical point or includes a connected component of $\partial\manif$. 
By Lemma~\ref{lm:cond_crlev_cont} the latter assumption implies that the mapping $\crvmap$ corresponding to each $\gfunc\in\OrbfMP$ its ordered set of exceptional values is continuous.

\subsection{Case $\RRR$.}
Notice that we have the following commutative diagram:
$$
\begin{CD}
\StabfM       @>{}>> \DiffM   @>{}>>      \DiffM/\StabfM        @>{}>>      \OrbfM  \\
 @V{}VV              @V{}VV                @V{}VV                            @V{}VV \\
\StabfMR      @>{}>> \DiffMRf  @>{}>>      \DiffMRf/\StabfMR      @>{}>>      \OrbfMR \\
@V{\prMPP}VV         @V{}VV                @V{}VV                    @V{\crvmap}VV \\
\DiffRfe  @>{}>> \DiffRf  @>{\padj}>> \DiffRf/\DiffRfe  @>{\adjd}>> \dntwo
\end{CD}
$$
in which $\adjd$ is a homeomorphism and the other right horizontal arrows are continuous bijections. Moreover, every upper vertical arrow is an embedding and every lower vertical one is a mapping onto.

The following lemma implies Theorem~\ref{th:Orbit} for the case $\Psp=\RRR$. 
\begin{lem}\label{lm:3_cond_on_f}
Suppose $\mrsfunc\in\smr$ satisfies $\condBndCrVal$.
Suppose also that 
\begin{enumerate}
\item the mapping $\crvmap:\OrbfMR\to\dntwo$ is continuous and there exist

\item a section $\sct:\dntwo\to\DiffRf$ of $\adjd\circ\padj$
(see Lemma~\ref{clm:constr_section_not_onto}), and

\item a section of $\prMPP$, i.e.\!\! a homomorphism $\difShift:\DiffRfe\to\DiffM$ such that $\difR\circ\mrsfunc=\mrsfunc\circ\difShift(\difR)$
(see Lemma~\ref{lm:pStMP_in_Pcrf}).

\end{enumerate}
Then there is a homeomorphism $\OrbfMR \approx \OrbfM \times \dntwo$.
\end{lem}
\proof
Consider the following composition:
$$
\begin{CD}
\sct\crvmap = \sct\circ\crvmap: \OrbfMR @>{\crvmap}>> \dntwo @>{\sct}>> \DiffRf.
\end{CD}
$$
Then each $\gfunc\in\OrbfMR$ can be represented in the form
$$
 \gfunc =\sct\crvmap(\gfunc)\ \circ \
 \bigl( \sct\crvmap(\gfunc)^{-1}   \circ \gfunc \bigr).
$$

\begin{claim}
$\sct\crvmap(\gfunc)^{-1} \circ \gfunc \in\OrbfM$.
\end{claim}
\proof
Since $\gfunc\in\OrbfMR$, it can be represented in the form
$\gfunc=\difR\circ\mrsfunc\circ\difM^{-1}$, where $\difR\in\DiffRf$ and $\difM\in\DiffM$.

Notice that $\sct\crvmap(\gfunc)^{-1} \circ \difR \in \DiffRfe$, whence by Proposition~\ref{pr:proof_A}
$$
\sct\crvmap(\gfunc)^{-1} \circ \gfunc  =
\bigl( \sct\crvmap(\gfunc)^{-1} \circ \difR \bigr) \circ \mrsfunc\circ\difM^{-1} =
\mrsfunc \circ
\difShift \left[\sct\crvmap(\gfunc)^{-1} \circ\difR \right] \circ \difM^{-1} \in \OrbfM. \qed
$$

Now it is easy to see that the following two mappings are continuous and mutually inverse:
$$
\begin{array}{lcl}
\homoodinv:\OrbfMR \to \OrbfM\times\dntwo, & \quad &
 \homoodinv(\gfunc) =
 \bigl(\,\sct\crvmap(\gfunc)^{-1} \circ \gfunc \,, \,\crvmap(\gfunc) \,\bigr), \\ [3mm]
\homood: \OrbfM\times\dntwo \to \OrbfMR, & \quad &
 \homood(\psi,z) = \sct(z)\circ \psi,
\end{array}
$$
where $\gfunc\in\OrbfMR$, $\psi\in\OrbfM$ and $z\in\dntwo$.
\endproof

\subsection{Case $\Psp=\Circle$, $\nv=0$.}
In this case $\mrsfunc$ is a locally trivial fibration over $\Circle$. 
Evidently, $\OrbfM\subset\OrbfMS$.
Conversely, suppose that  $\gfunc=\difR\circ\mrsfunc\circ\difM^{-1}\in\OrbfMS$ for $\difR\in\DiffSf$ and $\difM\in\DiffM$.
Then by Proposition~\ref{pr:proof_A_0} for $\difR$ there exists a diffeomorphism $\difM_1$ such that $\difR\circ\mrsfunc=\mrsfunc\circ\difM_1$, whence 
$\gfunc=\mrsfunc\circ\difM_1\circ\difM^{-1}\in\OrbfM$.
Thus $\OrbfM=\OrbfMS$.

\subsection{Case $\Psp=\Circle$, $\nv\geq1$.}
Define the mapping $\pfact:\DiffMSf\to\OrbfMS\times\impcr$ by
\[
\pfact(\difM,\difR) = \bigl( \difR\circ\mrsfunc\circ\difM^{-1}; 
\difR(1),\ldots,\difR(\nv) \bigr).
\]
Evidently, $\pfact$ is continuous.
Designate its image by $\pOrbfMS$: 
$$
 \pOrbfMS\,:=\,\pfact(\DiffMSf) \ \subset\ \OrbfMS\times\impcr.
$$ 

Let also \ $\prone:\pOrbfMS\to\OrbfMS$ \ be the restrictions of the standard projection $\OrbfMS\times\impcr\to\OrbfMS$ to $\pOrbfMS$.

Recall that $\ZZZ_{n}$ acts on $\impcr$ by cyclic shift $\genzn$ of coordinates:
$$ \genzn\cdot\{x_{\acnv}\}=\{x_{\acnv+1}\},$$
where for simplicity a point $(x_{1},\ldots,x_{\nv})\in\impcr$ is designated by $\{x_{\acnv}\}$.
This action together with the trivial action on $\OrbfMS$ yields an action of $\ZZZ_{n}$ on $\OrbfMS\times\impcr$.

By Theorem~\ref{th:Stab} $\prMPP(\StabfMS)/\DiffSfe$ is a cyclic group of some order $\kcnv$ dividing $\nv$.
Denote $\dcnv=\nv/\kcnv$ and define $\bdifR\in\DiffSf$ by 
\begin{equation}\label{equ:rotation_d}
\bdifR(\pnt)=\pnt+\dcnv\mod\nv.
\end{equation}
Let $\Zk$ be a subgroup of $\ZZZ_{\nv}$ generated by $\genzn^{\dcnv}$.

\begin{lem}\label{lm:struct_pOrbfMS}
$\pOrbfMS$ is invariant under $\genzn^{\dcnv}$ and $\prone$ yields a continuous bijection $\bij:\pOrbfMS/\Zk\to\OrbfMS$.
If $\crvmap$ is continuous, then $p_1$ is a $\kcnv$-sheet covering and $\bij$ is a homeomorphism.
\end{lem}
\proof
First we show that
$\bdifR\in\prMPP(\StabfMS)$, i.e.\! $(\bdifM,\bdifR)\in\StabfMS$ for some $\bdifM\in\DiffM$.

Indeed, since $\prMPP(\StabfMS)/\DiffSfe \approx \ZZZ_{\kcnv}$, there is 
$(\difM_1, \difR_1)\in\StabfMS$ such that $\difR_1(\acnv)=\acnv+\dcnv$ for all $\acnv=1,\ldots,\nv$.
Then $\bdifR^{-1}\circ\difR_1\in\DiffSfe$, whence 
$$
\mrsfunc = 
\bdifR\circ\bdifR^{-1}\circ\difR_1\circ\mrsfunc\circ\difM_1^{-1} = 
\bdifR\circ\mrsfunc\circ
\underbrace{\difShift(\bdifR^{-1}\circ\difR_1)\circ\difM_1^{-1}}_{\bdifM^{-1}}.
$$

Suppose that $(\gfunc, \{x_{\acnv}\})\in\pOrbfMS$, i.e.\!
$\gfunc=\difR\circ\mrsfunc\circ\difM^{-1}$ and $\difR(\acnv)=x_{\acnv}$ for $\acnv=1,\ldots,\nv$.
We have to show that $\genzn^{\dcnv}\cdot(\gfunc, \{x_{\acnv}\})=(\gfunc, \{x_{\acnv+\dcnv}\})$ also belongs to $\pOrbfMS$.
Then 
$$
\gfunc=\difR\circ\mrsfunc\circ\difM^{-1} = 
\difR\circ\bdifR\circ \mrsfunc\circ \bdifM^{-1}\circ\difM^{-1}.
$$

Moreover, 
$\difR\circ\bdifR(\acnv) =
\difR(\dcnv + \acnv)=x_{\dcnv+\acnv}$.
Thus putting $\hat\difR=\difR\circ\bdifR$ and 
$\hat\difM=\difM\circ\bdifM$ we see that 
$\pfact(\hat\difM,\hat\difR)=(\gfunc,\{x_{\dcnv + \acnv}\})$, 
i.e. $(\gfunc,\{x_{\dcnv + \acnv}\})$ belongs to the image $\pOrbfMS$ of $\pfact$.
Thus $\pOrbfMS$ is invariant under $\Zk$.
 
Since $\Zk$ is finite and its action is free and $p_1$-equivariant, it follows that the factor mapping $\pOrbfMS\to\pOrbfMS/\Zk$ is a $\kcnv$-sheet covering.
Moreover, since $p_1$ is onto, we obtain that $p_1$ yields a continuous bijection $\bij:\pOrbfMS/\Zk\to\OrbfMS$.

Finally, suppose that $\crvmap$ is continuous.
We will show that $\prone$ is a local homeomorphism. This will imply that so is $\bij$, whence $\bij$ is in fact a homeomorphism.
It suffices to prove that $\prone$ admits continuous local sections, i.e.\!
for every $(\gfunc,x)\in\pOrbfMS$ there exists a neighborhood $\Nbh_{\gfunc}$ of $\gfunc$ in $\OrbfMS$ and a continuous mapping $\Gfunc:\Nbh_{\gfunc}\to\impcr$ such that $(\gfunc',\Gfunc(\gfunc'))\in\pOrbfMS$ and $\Gfunc(\gfunc)=x$.

Notice that we have the following maps:
$$
\impcr \xrightarrow{\ \factn \ } \impcr/\ZZZ_{n} \xleftarrow{\ \crvmap \ } \OrbfMS.
$$

Let $(\gfunc,x)\in\pOrbfMS$ and $[x]=\factn(x)$ be the class of $x$ in $\impcr/\ZZZ_{n}$. 
Then $\crvmap(\gfunc)=[x]$.
Since $\factn$ is a covering, there is a neighborhood $\nbh_{x}$ of $x$ in $\impcr$ and a neighborhood $\anbh_{[x]}$ of $[x]$ in $\impcr/\ZZZ_{n}$ such that $\factn$ homeomorphically maps $\nbh_{x}$ onto $\anbh_{[x]}$.
Let $\Nbh_{\gfunc}=\crvmap^{-1}(\anbh_{[x]})$. 
Since $\crvmap$ is continuous, we obtain that $\Nbh_{\gfunc}$ is an open neighborhood of $\gfunc$. Then the mapping 
$\factn^{-1}\circ\crvmap:\Nbh_{\gfunc}\to\nbh_{x}$ is a local section of $\prone$.
\endproof

Consider now the projection \
$\crvmaptwo:\pOrbfMS\subset\OrbfMS\times\impcr\to\impcr.$ \
Since the composition 
$$
\crvmaptwo\circ\pfact:\DiffMSf \xrightarrow{\pfact} \OrbfMS\times\impcr \xrightarrow{\crvmaptwo} \impcr.
$$
is given by $(\difM,\difR)\mapsto\bigl(\difR(1),\ldots,\difR(\nv-1)\bigr)$, we get the following commutative diagram:
$$
\begin{CD}
\StabfM       @>{}>> \DiffM   @>{}>>      \DiffM/\StabfM        @>{}>>      \OrbfM  \\
 @V{}VV              @V{}VV                @V{}VV                            @V{}VV \\
\pStabfMS      @>{}>> \DiffMSf  @>{}>>      \DiffMSf/\pStabfMS      @>{}>>      \pOrbfMS \\ @V{\prMPP}VV         @V{}VV                @V{}VV                    @V{\crvmaptwo}VV \\
\DiffSfe  @>{}>> \DiffSf  @>{\padj}>> \DiffSf/\DiffSfe  @>{\adjd}>> \impcr
\end{CD}
$$
in which $\adjd$ is a homeomorphism and the other right horizontal arrows are continuous bijections, upper vertical and left horizontal arrows are embeddings, and lower vertical ones are surjective mappings.

Notice that 
\begin{enumerate}
\item 
$\crvmaptwo$ is continuous,
\item
$\adjd\circ\padj$ admit a continuous section $\sct$ (Lemma~\ref{lm:pStMP_in_Pcrf}), and
\item
$\prMPP$ admit a continuous section $\sectMPP$ being a homomorphism (Theorem~\ref{th:Stab}).
\end{enumerate}
Then by the arguments similar to the proof of Lemma~\ref{lm:3_cond_on_f}, it follows that the embedding $\OrbfM\equiv\OrbfM\times(1,\ldots,\nv)\subset\pOrbfMS$ extends to a homeomorphism
\[
\homood:\OrbfM\times\impcr \approx \pOrbfMS, \qquad
\homood(\gfunc,x)=(\sct(x)\circ\gfunc,x),
\]
where $\gfunc\in\OrbfM$ and $x\in\impcr$.

Finally, by Lemma~\ref{lm:struct_pOrbfMS} we have that $\OrbfMS=\pOrbfMS/\Zk$.
Since $\Zk$ trivially acts on $\OrbfM$, we obtain a homeomorphism:
\[
\OrbfMS \approx \OrbfM\times(\impcr/\Zk).
\]
Then it remains to apply Lemma~\ref{lm:FnZk_S1_symplex}.
Theorem~\ref{th:Orbit} is completed.
\qed
\section{Acknowledgements}
I am sincerely grateful to V.~V.~Sharko, D.~Bolotov, M.~Pankov,
E.~Polulyah, A.~Prishlyak, I.~Vlasenko for useful discussions.
I am indebted to K.~Feldman for careful reading the manuscript and for valuable comments.

\end{document}